\documentclass[12pt,english]{article}
\usepackage[T1]{fontenc}
\usepackage[latin9]{inputenc}
\usepackage[active]{srcltx}
\usepackage{babel}
\usepackage{amsmath}
\usepackage{amssymb}
\usepackage{wasysym}
\usepackage[bookmarks=false,
 breaklinks=false,pdfborder={0 0 1},backref=section,colorlinks=false]
 {hyperref}

\makeatletter

\providecommand{\tabularnewline}{\\}

\usepackage{wasysym}

\providecommand{\tabularnewline}{\\}

\usepackage{color}
\usepackage{graphicx}
\usepackage{amsfonts}
\usepackage{a4wide}
\usepackage{babel}
\usepackage{mathrsfs}

\providecommand{\U}[1]{\protect\rule{.1in}{.1in}}

\newtheorem{theorem}{Theorem}\newtheorem{corollary}[theorem]{Corollary}\newtheorem{definition}[theorem]{Definition}\newtheorem{example}[theorem]{Example}\newtheorem{lemma}[theorem]{Lemma}\newtheorem{proposition}[theorem]{Proposition}\newtheorem{remark}[theorem]{Remark}\newenvironment{proof}[1][Proof]{\noindent\textbf{#1.} }{\ \rule{0.5em}{0.5em}}

\usepackage[usenames,dvipsnames]{pstricks}
\usepackage{epsfig}
\usepackage{pst-grad}
\usepackage{pst-plot}

\usepackage{graphicx}
\usepackage{upgreek}

\usepackage{tikz}
\usepackage{hyphenat}

\makeatother

\begin{document}
\title{General monotonicity}
\author{M.D. Voisei}
\date{{}}

\maketitle
 
\begin{abstract}
This article employs techniques from convex analysis to present characterizations
of (maximal) $n-$monotonicity, similar to the well-established characterizations
of (maximal) monotonicity found in the existing literature. These
characterizations are further illustrated through examples. 
\end{abstract}
\textbf{Keywords} dual system, maximal cyclically monotone operator,
Fitzpatrick func\-tion, Rockafellar's antiderivative, planar rotation,
skew double-cone

\strut

\noindent\textbf{Mathematics Subject Classification (2020)} 47H05,
46N99, 47N10.

\section{Introduction and preliminaries }

The primary objective of this paper is to offer methods for identifying
(maximal) $n-$mon\-o\-tone operators via convex analysis techniques,
following the framework outlined in \cite{MR1009594,MR2207807,MR0193549},
where $n\in\mathbb{N}\cup\{+\infty)$, $n\ge1$, denoted as $n\in\overline{1,\infty}$.

In specific contexts, such as when the spaces are reflexive or when
the domain or range is convex, maximal cyclical monotonicity coincides
with maximal monotonicity (see e.g. \cite{MR1848999}). There is a
substantial body of literature addressing maximal monotonicity in
reflexive spaces. Consequently, our focus is on the general context
of dual systems that may not necessarily stem from reflexive spaces
or dualities between a space and its topological dual. The central
contribution of this article is the functional characterization of
maximal $n-$monotonicity within this broad context.

This article addresses topics that have been explored in several previous
works, including \cite{MR0275240,MR2286629,MR2291550,MR2369308}.
What distinguishes our approach is not only the extensive scope of
the context and the range of results presented, but also the introduction
of new, significant functions and operators that are deeply connected
to the theory of maximal $n-$monotone operators. These innovations
offer fresh insights and strengthen the link between convex analysis
and $n-$monotonicity.

In what follows, we adopt the following concepts and notations: for
a multi-valued operator (or multi-function) $T:X\rightrightarrows Y$
we associate it with its \emph{values}: $T(x)\subset Y$, $x\in X$,
\emph{graph}: $\operatorname*{Graph}T=\{(x,y)\in X\times Y\mid y\in T(x)\}$,
\emph{inverse:} $T^{-1}:Y\rightrightarrows X$, $\operatorname*{Graph}(T^{-1})=\{(y,x)\mid(x,y)\in\operatorname*{Graph}T\}$,
\emph{domain}: $D(T):=\{x\in X\mid T(x)\neq\emptyset\}=\Pr\nolimits_{X}(\operatorname*{Graph}T)$,
\emph{range}: $R(T):=\{y\in Y\mid y\in T(x),\ \mathrm{for\ some}\ x\in X\}=\Pr\nolimits_{Y}(\operatorname*{Graph}T)$,
\emph{direct--image: }for $A\subset X$, $T(A)=\cup_{x\in A}T(x)\subset Y$,
\emph{inverse--image:} for $B\subset Y$, $T^{-1}(B)=\{x\in X\mid T(x)\in B\}\subset X$.
Here $\Pr_{X}(x,y):=x$ and $\Pr_{Y}(x,y):=y$, $(x,y)\in X\times Y$
are the projections of $X\times Y$ onto $X$ and $Y$, respectively.

We do not identify an operator with its graph. Generally, a subset
$S\subset X\times Y$ is said to have a certain property if the operator
$T:X\rightrightarrows Y$, whose $\operatorname*{Graph}T=S$, has
that same property. Conversely, $T:X\rightrightarrows Y$ is said
to have a certain property if $\operatorname*{Graph}T\subset X\times Y$
possesses that property.

For a topologically vector space $(E,\mu)$, $A\subset E$, and $f,g:E\rightarrow\overline{\mathbb{R}}$
we denote by:

\medskip{}

\noindent$E^{*}$ -- the topological dual of $E$.

\medskip{}

\noindent$\operatorname*{cl}_{\mu}A$ -- the $\mu-$\emph{closure}
of $A$; $x\in\operatorname*{cl}_{\mu}A$ iff $x_{i}\stackrel{\mu}{\to}x$,
for some net $(x_{i})_{i}\subset A$. Here ``$\stackrel{\mu}{\to}$''
denotes the convergence of nets in the $\mu$ topology.

\noindent$\operatorname*{conv}A$ -- the \emph{convex hull} of $A$;
$x\in\operatorname*{conv}A$ iff $(x,1)={\displaystyle \sum_{k=1}^{n}}\lambda_{k}(x_{k},1)$
for some $(\lambda_{k})_{k\in\overline{1,n}}\subset\mathbb{R}^{+}:=[0,+\infty):=\{x\in\mathbb{R}\mid x\ge0\}$,
$(x_{k})_{k\in\overline{1,n}}\subset A$.

\medskip{}

\noindent$\operatorname*{Epi}f$ $:=\{(x,t)\in X\times\mathbb{R}\mid f(x)\leq t\}\subset X\times\mathbb{R}$
-- the \emph{epigraph} of $f$.

\medskip{}

\noindent$\operatorname*{cl}_{\mu}f$ -- the $\mu-$\emph{lower
semicontinuous hull} of $f$, which is the greatest $\mu-$lower semicontinuous
function majorized by $f$; $\operatorname*{Epi}(\operatorname*{cl}_{\mu}f)=\operatorname*{cl}_{\mu\times\tau_{0}}(\operatorname*{Epi}f)$,
where $\tau_{0}$ is the usual topology of $\mathbb{R}$.

\medskip{}

\noindent$\operatorname*{conv}f$ -- the \emph{convex hull} of $f$,
i.e., the greatest convex function majorized by $f$; \\
 $(\operatorname*{conv}f)(x):=\inf\{t\in\mathbb{R}\mid(x,t)\in\operatorname*{conv}(\operatorname*{Epi}f)\}$,
for $x\in X$; $\operatorname*{Epi}(\operatorname*{conv}f)=\operatorname*{conv}(\operatorname*{Epi}f)$.

\medskip{}

\noindent$\operatorname*{cl}_{\mu}\operatorname*{conv}f$ -- the
$\mu-$\emph{lower semicontinuous convex hull} of $f$, which is the
greatest $\mu-$lower semicontinuous convex function majorized by
$f$; $\operatorname*{Epi}(\operatorname*{cl}_{\mu}\operatorname*{conv}f):=\operatorname*{cl}_{\mu\times\tau_{0}}\operatorname*{conv}(\operatorname*{Epi}f)$.

\medskip{}

\noindent The set $[f\le g]$ $:=\{x\in E\mid f(x)\leq g(x)\}$; the
sets $[f=g]$, $[f<g]$, $[f>g]$, $[f\ge g]$ are similarly defined;
$f\le(<,>,\ge)g$ means $[f\le(<,>,\ge)g]=E$, or, equivalently, for
every $x\in E$, $f(x)\le(<,>,\ge)g(x)$.

\medskip{}

\noindent$\Lambda(E)$ -- the class of proper convex functions $f:E\rightarrow\overline{\mathbb{R}}$.
Recall that $f$ is \emph{proper} if $\operatorname*{dom}f:=\{x\in E\mid f(x)<\infty\}$
is nonempty and $f$ does not take the value $-\infty$; $f$ is \emph{convex}
if, for every $x,y\in E$, $t\in[0,1]$, we have $f(tx+(1-t)y))\le tf(x)+(1-t)f(y)$
while observing the conventions: $\infty-\infty=\infty$, $0\cdot\pm\infty=0$.

\medskip{}

\noindent$\Gamma_{\mu}(E)$ -- the class of functions $f\in\Lambda(E)$
that are $\mu-$lower semicontinuous.

\medskip{}

We avoid using the $\mu-$notation when the topology is implicitly
understood.

\medskip{}

A \emph{dual system} is a triple $(X,Y,c)$, where $X$, $Y$ are
vector spaces, and $c=\langle\cdot,\cdot\rangle:X\times Y\to\mathbb{R}$
is a bilinear map. The weakest topology on $X$, for which all linear
forms $\{X\ni x\to\langle x,y\rangle\mid y\in Y\}$ are continuous,
is denoted by $\sigma(X,Y)$ and it is called \emph{the weak topology
of} $X$ with respect to the duality $(X,Y,c)$. The weak topology
$\sigma(Y,X)$ on $Y$ is defined in a similar manner. The spaces
$X$ and $Y$ are regarded as locally convex spaces, endowed with
their respective weak topologies.

\eject

Given a dual system $(X,Y,\langle\cdot,\cdot\rangle)$, $S\subset X$,
and $f:X\rightarrow\overline{\mathbb{R}}$ 
\begin{itemize}
\item $f^{\ast}:Y\rightarrow\overline{\mathbb{R}}$ is the \emph{convex
conjugate} of $f:X\rightarrow\overline{\mathbb{R}}$ with respect
to the dual system $(X,Y,\langle\cdot,\cdot\rangle)$, $f^{\ast}(y):=\sup\{\left\langle x,y\right\rangle -f(x)\mid x\in X\}$
for $y\in Y$; similarly, the convex conjugate of $g:Y\rightarrow\overline{\mathbb{R}}$
is $g^{\ast}(x):=\sup\{\left\langle x,y\right\rangle -g(y)\mid y\in Y\}$
for $x\in X$; \vspace{0.5mm}
 
\item $\partial f(x)$ is the \emph{subdifferential} of $f:X\rightarrow\overline{\mathbb{R}}$
at $x\in X$; 
\begin{equation}
\partial f(x):=\{y\in Y\mid\forall w\in X,\ \left\langle w-x,y\right\rangle +f(x)\leq f(w)\},\label{eq:-11}
\end{equation}
for $x\in X$ with $f(x)\in\mathbb{R}$; $\partial f(x):=\emptyset$
whenever $f(x)\not\in\mathbb{R}$. Note that for an improper $f$
we have $\operatorname*{Graph}\partial f=\emptyset$. 
\item $Z:=X\times Y$ forms naturally a dual system $(Z,Z,\cdot)$ where
\begin{equation}
z\cdot w:=\langle x,v\rangle+\langle u,y\rangle,\ z=(x,y)\in Z,\ w=(u,v)\in Z.\label{eq:-161}
\end{equation}
\item $S^{\perp}:=\{y\in Y\mid\forall x\in S,\ \langle x,y\rangle=0\}$
is the \emph{orthogonal} \emph{of} $S$, $\sigma_{S}(y):=\iota_{S}^{*}(y)=\sup_{x\in S}\langle x,y\rangle$,
$y\in Y$ is the \emph{support functional }of $S$; note that $\sigma_{S}=\iota_{S^{\perp}}$
whenever $S$ is a linear subspace of $X$. 
\end{itemize}
The space $Z$ is endowed with a locally convex topology $\mu$ compatible
with the natural duality $(Z,Z,\cdot)$, meaning that $(Z,\mu)^{*}=Z$.
For example $\mu=\sigma(Z,Z)=\sigma(X,Y)\times\sigma(Y,X)$. All notions
associated with a proper function $f:Z\rightarrow\overline{\mathbb{R}}$
are similarly defined as above. Furthermore, with respect to the natural
dual system $(Z,Z,\cdot)$, the conjugate of $f$ is 
\begin{equation}
f^{\square}:Z\rightarrow\overline{\mathbb{R}},\quad f^{\square}(z)=\sup\{z\cdot z^{\prime}-f(z^{\prime})\mid z^{\prime}\in Z\}.\label{eq:-160}
\end{equation}
By the biconjugate formula, $f^{\square\square}=\operatorname*{cl}\operatorname*{conv}f$,
whenever $f^{\square}$ (or $\operatorname*{cl}\operatorname*{conv}f$)
is proper (see e.g. \cite{MR1921556}).

To $A\subset Z$ we associate $c_{A}:Z\rightarrow\overline{\mathbb{R}}$,
$c_{A}:=c+\iota_{A}$, where $\iota_{A}(z)=0$, for $z\in A$, $\iota_{A}(z)=+\infty$,
otherwise. The \emph{Fitzpatrick function} of $A$ (see \cite{MR1009594})
is 
\[
\varphi_{A}:Z\rightarrow\overline{\mathbb{R}},\ \varphi_{A}(z):=c_{A}^{\square}(z)=\sup\{z\cdot w-c(w)\mid w\in A\},
\]
and $\psi_{A}:Z\rightarrow\overline{\mathbb{R}}$, $\psi_{A}:=\varphi_{A}^{\square}=c_{A}^{\square\square}$.

In particular when $\varphi_{A}$ is proper (or $\operatorname*{cl}\operatorname*{conv}c_{A}$
is proper), 
\[
\psi_{A}=\operatorname*{cl}\operatorname*{conv}c_{A},\ \varphi_{A}=\psi_{A}^{\square}.
\]

Similarly, for a multi\-function $T:X\rightrightarrows Y$ we define
$c_{T}:=c_{\operatorname*{Graph}T}$, $\psi_{T}:=\psi_{\operatorname*{Graph}T}$,
and the \emph{Fitzpatrick function} of $T$, $\varphi_{T}:=\varphi_{\operatorname*{Graph}T}$.
By convention $\varphi_{\emptyset}=-\infty$, $c_{\emptyset}=\operatorname*{conv}c_{\emptyset}=\psi_{\emptyset}=+\infty$
in agreement with the usual conventions of $\inf_{\emptyset}=+\infty$,
$\sup_{\emptyset}=-\infty$.

In expanded form, for $T:X\rightrightarrows Y$ and $z=(x,y)\in Z$,
\begin{align}
\varphi_{T}(z) & =\sup\{z\cdot w-c(w)\mid w\in\operatorname*{Graph}T\}\nonumber \\
\varphi_{T}(x,y) & =\sup\{\langle x-u,v\rangle+\langle u,y\rangle\mid(u,v)\in\operatorname*{Graph}T\}.\label{def-Ff}
\end{align}

\eject

Let $(X,Y,\langle\cdot,\cdot\rangle)$ be a dual system and let $n\in\mathbb{N}$,
$n\ge1$. An operator $T:X\rightrightarrows Y$ is\emph{ } 
\begin{itemize}
\item \emph{$n-$monotone} (in $(X,Y,\langle\cdot,\cdot\rangle)$) if, for
every $\{(x_{i},y_{i})\}_{i\in\overline{1,n}}\subset\operatorname*{Graph}T$
with $x_{n+1}=x_{1}$ 
\begin{equation}
\sum_{i=1}^{n}\langle x_{i}-x_{i+1},y_{i}\rangle=\langle x_{1}-x_{2},y_{1}\rangle+\ldots+\langle x_{n-1}-x_{n},y_{n-1}\rangle+\langle x_{n}-x_{1},y_{n}\rangle\ge0.\label{eq:}
\end{equation}
Equivalently, using (\ref{eq:}) with the pairs in reversed order,
$T$ is $n-$monotone iff for every $\{(x_{i},y_{i})\}_{i\in\overline{1,n}}\subset\operatorname*{Graph}T$
with $x_{0}=x_{n}$ 
\begin{equation}
\begin{aligned}\langle x_{n}- & x_{n-1},y_{n}\rangle+\ldots+\langle x_{2}-x_{1},y_{2}\rangle+\langle x_{1}-x_{n},y_{1}\rangle\\
 & =\sum_{k=1}^{n}\langle x_{n+1-k}-x_{n-k},y_{n+1-k}\rangle=\sum_{i=1}^{n}\langle x_{i}-x_{i-1},y_{i}\rangle\ge0.
\end{aligned}
\label{eq:-8}
\end{equation}
\item \emph{cyclically or ($\infty-$)monotone} if, for every $k\in\mathbb{N}$,
$k\ge2$, $T$ is $k-$monotone. 
\item \emph{maximal $n-$monotone} if $T$ is $n-$monotone and it does
not allow for a proper $n-$mono\-tone extension in terms of graph
inclusion. Here $n\in\overline{1,\infty}$. 
\end{itemize}
Every operator is inherently $1-$monotone and $Z=X\times Y$ is the
only maximal $1-$monotone set. The term \textquotedbl monotone\textquotedbl{}
typically refers to $2-$mono\-tone operators. For $n,m\in\overline{2,\infty}$,
with $n\le m$, every $m-$monotone operator is also $n-$monotone.
Moreover, an operator that is both $m-$monotone and maximal $n-$monotone
is also maximal $m-$monotone. The prototype for cyclically monotone
operators is given by the convex subdifferential, as described in
\cite[Theorem 1, p. 500]{MR0193549} or Theorem \ref{RB} below.

\begin{definition} \emph{For $n\in\overline{1,\infty}$, the }$n-$th
Fitzpatrick function\emph{ associated with $T:X\rightrightarrows Y$
is denoted as $\overset{{\scriptscriptstyle \langle n\rangle}}{\varphi}_{T}:X\times Y\to\overline{\mathbb{R}}$
and is defined as follows: 
\begin{equation}
\overset{{\scriptscriptstyle \langle1\rangle}}{\varphi}_{T}(x,y):=\varphi_{T}(x,y):=\sup\{\langle x-u,v\rangle+\langle u,y\rangle\mid(u,v)\in\operatorname*{Graph}T\};\label{eq:-1}
\end{equation}
for $2\le n<\infty$, 
\begin{equation}
\begin{aligned}\overset{{\scriptscriptstyle \langle n\rangle}}{\varphi}_{T}(x,y): & =\negthickspace\negthickspace\negthickspace\sup_{\{(x_{i},y_{i})\}_{i\in\overline{1,n}}\subset\operatorname*{Graph}T}\negthickspace\negthickspace\negthickspace\negthickspace\negthickspace\negthickspace\negthickspace\negthickspace\negthickspace\{\langle x-x_{1},y_{1}\rangle+\langle x_{1}-x_{2},y_{2}\rangle+\ldots+\langle x_{n-1}-x_{n},y_{n}\rangle+\langle x_{n},y\rangle\};\\
 & =\sup_{\{(x_{i},y_{i})\}_{i\in\overline{1,n}}\subset\operatorname*{Graph}T}\{\langle x-x_{1},y_{1}\rangle+\sum_{i=2}^{n}\langle x_{i-1}-x_{i},y_{i}\rangle+\langle x_{n},y\rangle\};
\end{aligned}
\label{eq:-2}
\end{equation}
and 
\begin{equation}
\overset{{\scriptscriptstyle \langle\infty\rangle}}{\varphi}_{T}=\sup_{n\ge1}\overset{{\scriptscriptstyle \langle n\rangle}}{\varphi}_{T}.\label{eq:-13}
\end{equation}
}\end{definition}

Note that $\overset{{\scriptscriptstyle \langle1\rangle}}{\varphi}_{T}=c_{T}^{\square}$.
For $n\in\mathbb{N}$, $n\ge2$, after using a reversed pair order,
\begin{equation}
\begin{aligned}\overset{{\scriptscriptstyle \langle n\rangle}}{\varphi}_{T}(x,y):=\negthickspace\negthickspace\negthickspace & \sup_{\{(x_{k},y_{k})\}_{k\in\overline{1,n}}\subset\operatorname*{Graph}T}\negthickspace\negthickspace\negthickspace\negthickspace\negthickspace\negthickspace\{\langle x-x_{n},y_{n}\rangle+\sum_{k=2}^{n}\langle x_{n+2-k}-x_{n+1-k},y_{n+1-k}\rangle+\langle x_{1},y\rangle\}\\
= & \sup_{\{(x_{i},y_{i})\}_{i\in\overline{1,n}}\subset\operatorname*{Graph}T}\{\langle x-x_{n},y_{n}\rangle+\sum_{i=1}^{n-1}\langle x_{i+1}-x_{i},y_{i}\rangle+\langle x_{1},y\rangle\}.
\end{aligned}
\label{eq:-3}
\end{equation}
For $n,m\in\overline{1,\infty}$, with $n\le m$, it follows that
$\overset{{\scriptscriptstyle \langle n\rangle}}{\varphi}_{T}\le\overset{{\scriptscriptstyle \langle m\rangle}}{\varphi}_{T}$.
As a result $\overset{{\scriptscriptstyle \langle\infty\rangle}}{\varphi}_{T}=\lim_{n\to\infty}\overset{{\scriptscriptstyle \langle n\rangle}}{\varphi}_{T}$.

In general, for $n\in\mathbb{N}$, $n\ge1$, and $x\in D(T)$, $y\in Y$,
if we pick $(x_{i},y_{i})=(x,v)$, $k\in\overline{1,n}$, for some
$v\in T(x)$, then, from (\ref{eq:-1}), (\ref{eq:-2}), $\overset{{\scriptscriptstyle \langle n\rangle}}{\varphi}_{T}(x,y)\ge\langle x,y\rangle$.
Therefore 
\begin{equation}
\forall T:X\rightrightarrows Y,\forall n\in\overline{1,\infty},\ \operatorname*{Graph}T\!\subset(D(T)\times Y)\cup(X\times R(T))\subset[\overset{{\scriptscriptstyle \langle n\rangle}}{\varphi}_{T}\ge c]\!\subset\![\overset{{\scriptscriptstyle \langle\infty\rangle}}{\varphi}_{T}\ge c].\label{eq:-4}
\end{equation}

Notice also that, for \emph{$n\in\overline{1,\infty}$,} 
\begin{equation}
\overset{{\scriptscriptstyle \langle n\rangle}}{\varphi}_{T}(x,y)=\overset{{\scriptscriptstyle \langle n\rangle}}{\varphi}_{T^{-1}}(y,x),\ x\in X,\ y\in Y.\label{eq:-41}
\end{equation}

\begin{definition} \emph{For $n\in\overline{1,\infty}$, $\overset{{\scriptscriptstyle \langle n\rangle}}{\psi}_{T}:X\times Y\to\overline{\mathbb{R}}$
is the convex conjugate of the $n-$th Fitzpatrick function, defined
as 
\begin{equation}
\overset{{\scriptscriptstyle \langle n\rangle}}{\psi}_{T}=(\overset{{\scriptscriptstyle \langle n\rangle}}{\varphi}_{T})^{\square},\label{eq:-35}
\end{equation}
relative to the natural dual system $(Z,Z,\cdot)$.} \end{definition}

\begin{definition}\emph{ For $n\in\overline{1,\infty}$, consider
$\overset{{\scriptstyle \langle n\rangle}}{\pisces}_{T}:X\times Y\to\overline{\mathbb{R}}$
given by } 
\begin{equation}
\overset{{\scriptstyle \langle1\rangle}}{\pisces}_{T}=c_{T},\ \overset{{\scriptstyle \langle2\rangle}}{\pisces}_{T}(x,y)=\inf\{\langle x_{1},y\rangle+\langle x-x_{1},y_{2}\rangle\mid x_{1}\in T^{-1}(y),\ y_{2}\in T(x)\},\label{eq:-5}
\end{equation}
\begin{equation}
\begin{aligned}\mathrm{for}\ 3\le n<\infty,\ \overset{{\scriptstyle \langle n\rangle}}{\pisces}_{T}(x,y)=\inf\{\langle x_{1},y\rangle+\sum_{i=2}^{n-1}\langle x_{i}-x_{i-1},y_{i}\rangle & +\langle x-x_{n-1},y_{n}\rangle\mid\\
x_{1}\in T^{-1}(y),\ y_{n}\in T(x),\ \{(x_{i},y_{i})\}_{i\in\overline{2,n-1}} & \subset\operatorname*{Graph}T\},
\end{aligned}
\label{eq:-6}
\end{equation}
\begin{equation}
\mathrm{and}\ \overset{{\scriptstyle \langle\infty\rangle}}{\pisces}_{T}(x,y)=\inf_{n\ge1}\overset{{\scriptstyle \langle n\rangle}}{\pisces}_{T}(x,y).\label{eq:-7}
\end{equation}
\end{definition}

Note that $\operatorname*{dom}\overset{{\scriptstyle \langle1\rangle}}{\pisces}_{T}=\operatorname*{Graph}T$
and, for $n\in\overline{2,\infty}$, $\operatorname*{dom}\overset{{\scriptstyle \langle n\rangle}}{\pisces}_{T}=[\overset{{\scriptstyle \langle n\rangle}}{\pisces}_{T}<\infty]=D(T)\times R(T)$.

\begin{lemma} \label{HT} Let $(X,Y,\langle\cdot,\cdot\rangle)$
be a dual system, let $n\in\overline{1,\infty}$, and let $T:X\rightrightarrows Y$.
Then 
\begin{equation}
\forall n\in\overline{1,\infty},\ \overset{{\scriptstyle \langle n\rangle}}{\varphi}_{T}=(\overset{{\scriptstyle \langle n\rangle}}{\pisces}_{T})^{\square},\ \overset{{\scriptstyle \langle n\rangle}}{\psi}_{T}=(\overset{{\scriptstyle \langle n\rangle}}{\pisces}_{T})^{\square\square}\label{eq:-12}
\end{equation}
and, whenever $\overset{{\scriptstyle \langle n\rangle}}{\varphi}_{T}$
or $\operatorname*{cl}\operatorname*{conv}(\overset{{\scriptstyle \langle n\rangle}}{\pisces}_{T})$
is proper, $\overset{{\scriptstyle \langle n\rangle}}{\psi}_{T}=\operatorname*{cl}\operatorname*{conv}(\overset{{\scriptstyle \langle n\rangle}}{\pisces}_{T})$.

For $n,m\in\overline{1,\infty}$, with $n\le m$, $\overset{{\scriptstyle \langle m\rangle}}{\pisces}_{T}\le\overset{{\scriptstyle \langle n\rangle}}{\pisces}_{T}$;
in particular 
\begin{equation}
\forall n\in\overline{1,\infty},\ \overset{{\scriptstyle \langle n\rangle}}{\psi}_{T}\le\overset{{\scriptstyle \langle n\rangle}}{\pisces}_{T}\le\overset{{\scriptstyle \langle1\rangle}}{\pisces}_{T}=c_{T},\ \operatorname*{Graph}T\subset[\overset{{\scriptstyle \langle n\rangle}}{\pisces}_{T}\le c]\subset[\overset{{\scriptstyle \langle n\rangle}}{\psi}_{T}\le c].\label{eq:-17}
\end{equation}

\end{lemma}

\begin{proof} For (\ref{eq:-12}) we use (\ref{eq:-1}), (\ref{eq:-2}),
(\ref{eq:-13}), (\ref{eq:-5}), (\ref{eq:-6}), (\ref{eq:-7}) followed
by the biconjugate formula.

For every $(x,y)\in\operatorname*{Graph}T$, pick $x_{1}=x$ or $y_{2}=y$
in the definition of $\overset{{\scriptstyle \langle2\rangle}}{\pisces}_{T}$
to get $\overset{{\scriptstyle \langle2\rangle}}{\pisces}_{T}\le c_{T}=\overset{{\scriptstyle \langle1\rangle}}{\pisces}_{T}$.

Let $n,m\in\mathbb{N}$ be such that $2\le n<m$. For every $(x,y)\in D(T)\times R(T)$
and every $x_{1}\in T^{-1}(y)$, $y_{m}\in T(x)$, $\{(x_{i},y_{i})\}_{i\in\overline{2,m-1}}\subset\operatorname*{Graph}T$,
\[
\text{\ensuremath{\overset{{\scriptstyle \langle m\rangle}}{\pisces}_{T}}}(x,y)\le\langle x_{1},y\rangle+\sum_{i=2}^{m-1}\langle x_{i}-x_{i-1},y_{i}\rangle+\langle x-x_{m-1},y_{m}\rangle.
\]
For $k\in\overline{n,m-1}$ take $(x_{k},y_{k})=(x_{n-1},y_{n-1})$
to get 
\[
\forall(x,y)\in D(T)\times R(T),\ \forall x_{1}\in T^{-1}(y),\ \forall y_{n}\in T(x),\ \forall\{(x_{i},y_{i})\}_{i\in\overline{2,n-1}}\subset\operatorname*{Graph}T,
\]
\[
\text{\ensuremath{\overset{{\scriptstyle \langle m\rangle}}{\pisces}_{T}}}(x,y)\le\langle x_{1},y\rangle+\sum_{i=2}^{n-1}\langle x_{i}-x_{i-1},y_{i}\rangle+\langle x-x_{n-1},y_{n}\rangle.
\]
Pass to infimum over $x_{1}\in T^{-1}(y)$, $y_{n}\in T(x)$, $\{(x_{i},y_{i})\}_{i\in\overline{2,n-1}}\subset\operatorname*{Graph}T$
to find 
\[
\forall(x,y)\in D(T)\times R(T),\ \text{\ensuremath{\overset{{\scriptstyle \langle m\rangle}}{\pisces}_{T}}}(x,y)\le\text{\ensuremath{\overset{{\scriptstyle \langle n\rangle}}{\pisces}_{T}}}(x,y).
\]
Since elsewhere $\text{\ensuremath{\overset{{\scriptstyle \langle m\rangle}}{\pisces}_{T}}}(x,y)=\text{\ensuremath{\overset{{\scriptstyle \langle n\rangle}}{\pisces}_{T}}}(x,y)=+\infty$
the proof is complete. \end{proof}

\begin{lemma} \label{rep} Let $(X,Y,\langle\cdot,\cdot\rangle)$
be a dual system, and let $f\in\Gamma(Z)$ be such that $f\ge c$.
Then 
\begin{equation}
[f=c]\subset[f^{\square}=c].\label{eq:-19}
\end{equation}
\end{lemma}

\begin{proof} Let $z\in[f=c]$. Then 
\[
\forall w\in Z,\ f'(z;w):=\lim_{t\downarrow0}\tfrac{1}{t}(f(z+tw)-f(z))\ge\lim_{t\downarrow0}\tfrac{1}{t}(c(z+tw)-c(z))=z\cdot w,
\]
which shows that $z\in\partial f(z)$, because $f$ is convex. Hence
$f(z)+f^{\square}(z)=z\cdot z=2c(z)$ and so, $z\in[f^{\square}=c]$.
\end{proof}

\eject

\section{$n-$Monotonicity}

\begin{theorem} \label{cm0} Let $(X,Y,\langle\cdot,\cdot\rangle)$
be a dual system and let $n\in\overline{1,\infty}$.

The operator $T:X\rightrightarrows Y$ is $n-$monotone iff $\overset{{\scriptstyle \langle n\rangle}}{\pisces}_{T}\ge c$.

In this case, $\operatorname*{Graph}T\subset[\overset{{\scriptstyle \langle n\rangle}}{\pisces}_{T}=c]$
and, whenever $\operatorname*{Graph}T\not=\emptyset$, $\overset{{\scriptstyle \langle n\rangle}}{\pisces}_{T}$
is proper. \end{theorem}

\begin{proof} If $n=1$ or $\operatorname*{Graph}T=\emptyset$, then
the conclusion is straightforward, as $\overset{{\scriptstyle \langle1\rangle}}{\pisces}_{T}=c_{T}$.

Otherwise, suppose that $n\in\overline{2,\infty}$ and $T:X\rightrightarrows Y$
has $\operatorname*{Graph}T\not=\emptyset$.

Assume first that $n$ is finite. If $T$ is $n-$monotone then, for
every $(x,y)\in D(T)\times R(T)$ and every $x_{1}\in T^{-1}(y),\ y_{n}\in T(x),\ \{(x_{i},y_{i})\}_{i\in\overline{2,n-1}}\subset\operatorname*{Graph}T$,
the $n-$monotonicity of $T$ on the pairs $(x_{1},y)$, $(x,y_{n})$,
$(x_{n-1},y_{n-1}),\ldots,(x_{2},y_{2})$ provides 
\[
\langle x_{1}-x,y\rangle+\langle x-x_{n-1},y_{n}\rangle+\sum_{i=2}^{n-1}\langle x_{i}-x_{i-1},y_{i}\rangle\ge0.
\]
Pass to infimum over $x_{1}\in T^{-1}(y),\ y_{n}\in T(x),\ \{(x_{i},y_{i})\}_{i\in\overline{2,n-1}}\subset\operatorname*{Graph}T$
to obtain $(\overset{{\scriptstyle \langle n\rangle}}{\pisces}_{T}-c)(x,y)\ge0$.
Therefore $\overset{{\scriptstyle \langle n\rangle}}{\pisces}_{T}\ge c$.

Conversely, suppose that $\overset{{\scriptstyle \langle n\rangle}}{\pisces}_{T}\ge c$
. For every $\{(x_{i},y_{i})\}_{i\in\overline{1,n}}\subset\operatorname*{Graph}T$,
\[
\sum_{i=2}^{n}\langle x_{i}-x_{i-1},y_{i}\rangle+\langle x_{1}-x_{n},y_{1}\rangle\ge(\overset{{\scriptstyle \langle n\rangle}}{\pisces}_{T}-c)(x_{n},y_{1})\ge0.
\]
According to (\ref{eq:-8}), $T$ is $n-$monotone.

The case $n=\infty$ follows directly as a consequence of the cases
where $1\le n<\infty$.

In general, for every $n\in\overline{1,\infty}$, $\overset{{\scriptstyle \langle n\rangle}}{\pisces}_{T}\le\overset{{\scriptstyle \langle1\rangle}}{\pisces}_{T}=c_{T}$
which implies that $\operatorname*{Graph}T\subset[\overset{{\scriptstyle \langle n\rangle}}{\pisces}_{T}\le c]$.
When $T$ is $n-$monotone, meaning that $\overset{{\scriptstyle \langle n\rangle}}{\pisces}_{T}\ge c$,
we obtain $\operatorname*{Graph}T\subset[\overset{{\scriptstyle \langle n\rangle}}{\pisces}_{T}=c]$.
\end{proof}

\begin{theorem} \label{cm1} Let $(X,Y,\langle\cdot,\cdot\rangle)$
be a dual system and let $n,m\in\overline{1,\infty}$. The operator
$T:X\rightrightarrows Y$ is $(n+m)-$monotone iff one of the following
conditions holds:

\medskip{}

\noindent\emph{(i) }$\overset{{\scriptstyle \langle n\rangle}}{\varphi}_{T}\le\overset{{\scriptstyle \langle m\rangle}}{\pisces}_{T}$,
$\qquad$ \emph{(ii)} $\overset{{\scriptstyle \langle n\rangle}}{\varphi}_{T}\le\operatorname*{conv}\overset{{\scriptstyle \langle m\rangle}}{\pisces}_{T}$,
$\qquad$ \emph{(iii)} $\overset{{\scriptstyle \langle n\rangle}}{\varphi}_{T}\le\overset{{\scriptstyle \langle m\rangle}}{\psi}_{T}$.

\medskip{}

In this case, 
\begin{equation}
\forall k\in\overline{1,n+m-1},\ \operatorname*{Graph}T\subset[\overset{{\scriptstyle \langle k\rangle}}{\varphi}_{T}=c]\cap[\overset{{\scriptstyle \langle k\rangle}}{\pisces}_{T}=c]\cap[\operatorname*{conv}\overset{{\scriptstyle \langle k\rangle}}{\pisces}_{T}=c]\cap[\overset{{\scriptstyle \langle k\rangle}}{\psi}_{T}=c],\label{eq:-36}
\end{equation}
and, whenever $\operatorname*{Graph}T\neq\emptyset$, $\overset{{\scriptstyle \langle k\rangle}}{\varphi}_{T}$,
$\overset{{\scriptstyle \langle k\rangle}}{\psi}_{T}$, $\operatorname*{conv}\overset{{\scriptstyle \langle k\rangle}}{\pisces}_{T}$,
$\overset{{\scriptstyle \langle k\rangle}}{\pisces}_{T}$ are proper,
$k\in\overline{1,n+m-1}$. \end{theorem}

\begin{proof} If $\operatorname*{Graph}T=\emptyset$, then $\overset{{\scriptstyle \langle n\rangle}}{\varphi}_{T}=-\infty$,
$\overset{{\scriptstyle \langle m\rangle}}{\pisces}_{T}=\operatorname*{conv}\overset{{\scriptstyle \langle m\rangle}}{\pisces}_{T}=\overset{{\scriptstyle \langle m\rangle}}{\psi}_{T}=+\infty$,
and the conclusion follows. Otherwise, suppose that $\operatorname*{Graph}T\not=\emptyset$.

Assume that $T$ is $(n+m)-$monotone, where $m\ge1$, $n\ge1$ are
finite integers.

For every $(x,y)\in D(T)\times R(T)$, $x_{1}'\in T^{-1}(y)$, $y_{m}'\in T(x)$
and every $\{(x_{i},y_{i})\}_{i\in\overline{1,n}}\cup\{(x'_{j},y'_{j})\}_{j\in\overline{2,m-1}}\subset\operatorname*{Graph}T$,
the $(n+m)-$monotonicity condition (\ref{eq:}) for the pairs 
\[
(x_{n},y_{n}),(x_{n-1},y_{n-1}),\ldots,(x_{1},y_{1}),\ \ (x,y'_{m}),(x'_{m-1},y'_{m-1}),\ldots,(x'_{2},y'_{2}),(x'_{1},y)
\]
provides 
\[
\sum_{i=2}^{n}\langle x_{i}-x_{i-1},y_{i}\rangle+\langle x_{1}-x,y_{1}\rangle+\langle x-x_{m-1}',y'_{m}\rangle+\sum_{j=2}^{m-1}\langle x'_{j}-x'_{j-1},y'_{j}\rangle+\langle x_{1}'-x_{n},y\rangle\ge0,
\]
\[
\langle x-x_{1},y_{1}\rangle+\sum_{i=2}^{n}\langle x_{i-1}-x_{i},y_{i}\rangle+\langle x_{n},y\rangle\le\langle x_{1}',y\rangle+\sum_{j=2}^{m-1}\langle x'_{j}-x'_{j-1},y'_{j}\rangle+\langle x-x_{m-1}',y'_{m}\rangle.
\]

Pass to supremum over $\{(x_{k},y_{k})\}_{k\in\overline{1,n}}\subset\operatorname*{Graph}T$
on the left hand side and to infimum over $x_{1}'\in T^{-1}(y)$,
$y_{m}'\in T(x)$, $\{(x'_{j},y'_{j})\}_{j\in\overline{2,m-1}}\subset\operatorname*{Graph}T$
on the right hand side to get $\overset{{\scriptstyle \langle n\rangle}}{\varphi}_{T}(x,y)\le\overset{{\scriptstyle \langle m\rangle}}{\pisces}_{T}(x,y)$.
Since $\overset{{\scriptstyle \langle m\rangle}}{\pisces}_{T}=+\infty$
outside $D(T)\times R(T)$, we infer that (i) is true. In particular,
$\operatorname*{Graph}T\subset[\overset{{\scriptstyle \langle n\rangle}}{\varphi}_{T}=c]\cap[\overset{{\scriptstyle \langle m\rangle}}{\pisces}_{T}=c]$,
since, according to (\ref{eq:-4}), (\ref{eq:-17}), $\operatorname*{Graph}T\subset[\overset{{\scriptstyle \langle n\rangle}}{\varphi}_{T}\ge c]\cap[\overset{{\scriptstyle \langle m\rangle}}{\pisces}_{T}\le c]$,
and $\overset{{\scriptstyle \langle n\rangle}}{\varphi}_{T}$ is proper.

Assume that (i) holds. For every $\{(x_{i},y_{i})\}_{i\in\overline{1,n}}\cup\{(x'_{j},y'_{j})\}_{j\in\overline{1,m}}\subset\operatorname*{Graph}T$
\[
\langle x_{m}'-x_{1},y_{1}\rangle+\sum_{i=2}^{n}\langle x_{i-1}-x_{i},y_{i}\rangle+\langle x_{n},y_{1}'\rangle\le\overset{{\scriptstyle \langle n\rangle}}{\varphi}_{T}(x_{m}',y_{1}')\le\overset{{\scriptstyle \langle m\rangle}}{\pisces}_{T}(x_{m}',y_{1}')
\]
\[
\le\langle x_{1}',y_{1}'\rangle+\sum_{j=2}^{m-1}\langle x'_{j}-x'_{j-1},y'_{j}\rangle+\langle x_{m}'-x_{m-1}',y'_{m}\rangle=\langle x_{1}',y_{1}'\rangle+\sum_{j=2}^{m}\langle x'_{j}-x'_{j-1},y'_{j}\rangle,
\]
which leads to 
\[
\sum_{i=2}^{n}\langle x_{i}-x_{i-1},y_{i}\rangle+\langle x_{1}-x'_{m},y_{1}\rangle+\sum_{j=2}^{m}\langle x'_{j}-x'_{j-1},y'_{j}\rangle+\langle x_{1}'-x_{n},y'_{1}\rangle\ge0,
\]
i.e., the $(n+m)-$monotonicity condition for the pairs 
\[
(x_{n},y_{n}),(x_{n-1},y_{n-1}),\ldots,(x_{1},y_{1}),\ \ (x'_{m},y'_{m}),(x'_{m-1},y'_{m-1}),\ldots,(x'_{1},y'_{1}),
\]
holds. Therefore $T$ is $(n+m)-$monotone.

(i) $\Leftrightarrow$ (ii) $\Leftrightarrow$ (iii) are plain, $\overset{{\scriptstyle \langle m\rangle}}{\psi}_{T}$,
$\operatorname*{conv}\overset{{\scriptstyle \langle m\rangle}}{\pisces}_{T}$,
$\overset{{\scriptstyle \langle m\rangle}}{\pisces}_{T}$ are proper,
and (\ref{eq:-36}) is true in this case, since $\overset{{\scriptstyle \langle n\rangle}}{\varphi}_{T}\in\Gamma(Z)$
and $\overset{{\scriptstyle \langle n\rangle}}{\varphi}_{T}\le\overset{{\scriptstyle \langle m\rangle}}{\psi}_{T}=\operatorname*{cl}\operatorname*{conv}(\overset{{\scriptstyle \langle m\rangle}}{\pisces}_{T})\le\operatorname*{conv}\overset{{\scriptstyle \langle m\rangle}}{\pisces}_{T}\le\overset{{\scriptstyle \langle m\rangle}}{\pisces}_{T}$.

In case $m=\infty$ (or $n=\infty$) we use the above argument for
any finite integer $1\le m<\infty$ ($1\le n<\infty$) and pass to
infimum (supremum) over $m\ge1$ ($n\ge1$) to conclude. \end{proof}

\begin{corollary} \label{cm2} Let $(X,Y,\langle\cdot,\cdot\rangle)$
be a dual system and let $n\in\overline{1,\infty}$. The operator
$T:X\rightrightarrows Y$ is $(n+1)-$monotone iff one of the following
conditions holds:

\medskip{}

\noindent\emph{(i)} $\overset{{\scriptscriptstyle \langle n\rangle}}{\varphi}_{T}\le c_{T}$,
$\qquad$ \emph{(ii)} $\operatorname*{Graph}T\subset[\overset{{\scriptscriptstyle \langle n\rangle}}{\varphi}_{T}\le c]$,
$\qquad$ \emph{(iii)} $\operatorname*{Graph}T\subset[\overset{{\scriptscriptstyle \langle n\rangle}}{\varphi}_{T}=c]$.
\end{corollary}

\noindent\begin{proof} We use Theorem \ref{cm1} for $m=1$ to infer
that $T$ is $(n+1)-$monotone iff $\overset{{\scriptscriptstyle \langle n\rangle}}{\varphi}_{T}\le\overset{{\scriptstyle \langle1\rangle}}{\pisces}_{T}=c_{T}$,
that is, (i) holds; (i) $\Leftrightarrow$ (ii) is plain; and (ii)
$\Leftrightarrow$ (iii) follows with the aid of (\ref{eq:-4}). \end{proof}

\begin{definition} \emph{Given $(X,Y,\langle\cdot,\cdot\rangle)$
a dual system, $n\in\mathbb{N}$, $n\ge2$, and $T:X\rightrightarrows Y$,
the pair $(x,y)\in X\times Y$ is $n-$}monotonically related\emph{
(m.r. for short) }to\emph{ $T$ if, for every $\{(x_{i},y_{i})\}_{i\in\overline{1,n-1}}\subset\operatorname*{Graph}T$,
\begin{equation}
\langle x-x_{n-1},y\rangle+\langle x_{n-1}-x_{n-2},y_{n-1}\rangle+\ldots+\langle x_{2}-x_{1},y_{2}\rangle+\langle x_{1}-x,y_{1}\rangle\ge0.\label{eq:-9}
\end{equation}
The pair $(x,y)\in X\times Y$ is cyclically $(\infty-)$}monotonically
related\emph{ (m.r. for short) }to\emph{ $T$ if, for every $n\ge2$,
$(x,y)$ is $n-$m.r. to $T$. }

\medskip{}

\emph{For $n\in\overline{2,\infty}$, if $T$ is $n-$monotone, then
$(x,y)$ is $n-$m.r. to $T$ iff $\operatorname*{Graph}T\cup\{(x,y)\}$
is $n-$monotone.} \end{definition}

\begin{theorem} \label{cm3} Let $(X,Y,\langle\cdot,\cdot\rangle)$
be a dual system, let $n\in\overline{2,\infty}$, and let $T:X\rightrightarrows Y$.
Then

\medskip{}

\emph{(i)} $[\overset{{\scriptscriptstyle \langle n-1\rangle}}{\varphi}_{T}\le c]=\{(x,y)\in X\times Y\mid(x,y)\ \mathrm{is}\ n-\mathrm{m.r.\ to}\ T\}$;

\medskip{}

\emph{(ii)} $T$ is $n-$monotone iff $\operatorname*{Graph}T\subset[\overset{{\scriptscriptstyle \langle n-1\rangle}}{\varphi}_{T}\le c]$
iff $\operatorname*{Graph}T\subset[\overset{{\scriptscriptstyle \langle n-1\rangle}}{\varphi}_{T}=c]$.
\end{theorem}

\noindent\begin{proof} (i) From (\ref{eq:-2}), (\ref{eq:-9}),
is easily checked that, for $n\in\mathbb{N}$, $n\ge2$, $(x,y)$
is $n-$m.r. to $T$ iff $(x,y)\in[\overset{{\scriptscriptstyle \langle n-1\rangle}}{\varphi}_{T}\le c]$.
As a consequence the same holds for $n=\infty$.

(ii) This fact has already been observed in Corollary \ref{cm2}.
Alternatively, it follows from (i) and (\ref{eq:-4}). \end{proof}

\strut

The concept of ($2-$)monotonicity has been extensively studied in
\cite{MR2207807,MR2389004,MR2453098,MR2577332,MR2583911,MR2594359,MR2899842,MR3917361}
with the key findings summarized in the following theorem.

\begin{theorem} \ref{m2}\label{m2} Let $(X,Y,\langle\cdot,\cdot\rangle)$
be a dual system and let $T:X\rightrightarrows Y$. The operator $T:X\rightrightarrows Y$
is ($2-$)monotone iff one of the following conditions holds:

\medskip{}

\begin{tabular}{ll}
\emph{(i)} $\operatorname*{Graph}T\subset[\varphi_{T}=c]$,  & \emph{(v)} $\operatorname*{conv}c_{T}\ge c$,\tabularnewline
\emph{(ii)} $\overset{{\scriptscriptstyle \langle1\rangle}}{\varphi}_{T}\le c_{T}$,  & \emph{(vi)} there is $h\in\Lambda(Z)$ such that $h\ge c$ and $\operatorname*{Graph}T\subset[h=c]$,\tabularnewline
\emph{(iii)} $\overset{{\scriptscriptstyle \langle1\rangle}}{\varphi}_{T}\le\overset{{\scriptscriptstyle \langle1\rangle}}{\psi}_{T}$,  & \emph{(vii)} $\overset{{\scriptscriptstyle \langle1\rangle}}{\varphi}_{T^{+}}\ge c$,\tabularnewline
\emph{(iv)} $\overset{{\scriptscriptstyle \langle1\rangle}}{\psi}_{T}\ge c$,  & \emph{(viii)} $T^{++}$ is monotone. \tabularnewline
\end{tabular}

\medskip{}

Here $\operatorname*{Graph}(T^{+})=[\overset{{\scriptscriptstyle \langle1\rangle}}{\varphi}_{T}\le c]$
and $T^{++}=(T^{+})^{+}$. \end{theorem}

\begin{theorem} \label{m-inf} Let $(X,Y,\langle\cdot,\cdot\rangle)$
be a dual system and let $T:X\rightrightarrows Y$. The operator $T:X\rightrightarrows Y$
is cyclically monotone iff one of the following conditions holds:

\begin{tabular}{ll}
\emph{(i)} $\operatorname*{Graph}T\subset[\overset{{\scriptscriptstyle \langle\infty\rangle}}{\varphi}_{T}(\le)=c]$,  & \emph{(iv)} $\overset{{\scriptscriptstyle \langle\infty\rangle}}{\psi}_{T}\ge c$, \tabularnewline
\emph{(ii)} $\overset{{\scriptscriptstyle \langle\infty\rangle}}{\varphi}_{T}\le\overset{{\scriptscriptstyle \langle\infty\rangle}}{\pisces}_{T}$,  & \emph{(v)} $\operatorname*{conv}\overset{{\scriptscriptstyle \langle\infty\rangle}}{\pisces}_{T}\ge c$,\tabularnewline
\emph{(iii)} $\overset{{\scriptscriptstyle \langle\infty\rangle}}{\varphi}_{T}\le\overset{{\scriptscriptstyle \langle\infty\rangle}}{\psi}_{T}$,  & \emph{(vi)} there is $h\in\Lambda(Z)$ such that $c\le h\le\overset{{\scriptscriptstyle \langle\infty\rangle}}{\pisces}_{T}$.\tabularnewline
\end{tabular}\end{theorem}

\noindent\begin{proof} We may assume that $\operatorname*{Graph}T\not=\emptyset$,
otherwise the argument is plain.

In this case and under any of the given conditions, $\overset{{\scriptscriptstyle \langle\infty\rangle}}{\varphi}_{T}$,
$\overset{{\scriptscriptstyle \langle\infty\rangle}}{\psi}_{T}$,
$\overset{{\scriptscriptstyle \langle\infty\rangle}}{\pisces}_{T}$
are proper.

\noindent The equivalencies of $T$ being cyclically monotone with
(i), (ii), (iii) are consequences of Theorems \ref{cm3}, \ref{cm1}
for $n=m=\infty$.

(iii) $\Rightarrow$ (iv) For every $z\in Z$, 
\[
\overset{{\scriptscriptstyle \langle\infty\rangle}}{\psi}_{T}(z)\ge\tfrac{1}{2}(\overset{{\scriptscriptstyle \langle\infty\rangle}}{\psi}_{T}(z)+\overset{{\scriptscriptstyle \langle\infty\rangle}}{\varphi}_{T}(z))=\tfrac{1}{2}((\overset{{\scriptscriptstyle \langle\infty\rangle}}{\varphi}_{T})^{\square}(z)+\overset{{\scriptscriptstyle \langle\infty\rangle}}{\varphi}_{T}(z))\ge\tfrac{1}{2}z\cdot z=c(z).
\]

(iv) $\Rightarrow$ (v) follows from $\overset{{\scriptscriptstyle \langle\infty\rangle}}{\psi}_{T}=\operatorname*{cl}\operatorname*{conv}\overset{{\scriptscriptstyle \langle\infty\rangle}}{\pisces}_{T}\le\operatorname*{conv}\overset{{\scriptscriptstyle \langle\infty\rangle}}{\pisces}_{T}$.

(v) $\Rightarrow$ (vi) is straightforward with $h=\operatorname*{conv}\overset{{\scriptscriptstyle \langle\infty\rangle}}{\pisces}_{T}$.

If (vi) holds then $T$ is cyclically monotone due to Theorem \ref{cm0}.
\end{proof}

\strut

Notice that condition (vi) in both Theorems \ref{m2}, \ref{m-inf}
shares the same structure, because $h\ge c$ and $\operatorname*{Graph}T\subset[h=c]$
iff $c\le h\le c_{T}=\overset{{\scriptscriptstyle \langle1\rangle}}{\pisces}_{T}$.

\begin{remark} For $n\in\overline{2,\infty}$, let $T_{n}^{+}:Z\rightrightarrows Z$
be defined by $\operatorname*{Graph}(T_{n}^{+})=[\overset{{\scriptscriptstyle \langle n-1\rangle}}{\varphi}_{T}\le c]$,
i.e., $z\in\operatorname*{Graph}(T_{n}^{+})$ iff $z$ is $n-$m.r.
to $T$, and let $T_{n}^{++}=(T_{n}^{+})_{n}^{+}$.

For any operator $T$, we have $\operatorname*{Graph}T\subset\operatorname*{Graph}(T_{(2)}^{++})$,
However, for $n\in\overline{3,\infty}$, it is generally not true
that $\operatorname*{Graph}T\subset\operatorname*{Graph}(T_{n}^{++})$.
More precisely, 
\begin{equation}
\exists T:\mathbb{R}^{2}\rightrightarrows\mathbb{R}^{2},\ \forall n\in\overline{3,\infty},\ \operatorname*{Graph}T\not\subset\operatorname*{Graph}(T_{n}^{++}).\label{eq:-14}
\end{equation}

Take $T:\mathbb{R}^{2}\rightrightarrows\mathbb{R}^{2}$ given by $\operatorname*{Graph}T=\{(0,0),(5,5)\}$.
Then $\{(1,1),(2,0.32)\}\subset\operatorname*{Graph}(T_{\infty}^{+})=\cap_{n=2}^{\infty}\operatorname*{Graph}(T_{n}^{+})$,
because 
\[
\operatorname*{Graph}T\cup\{(1,1)\}\subset\operatorname*{Graph}(\partial(\tfrac{1}{2}x^{2}))=\{(x,x)\mid x\in\mathbb{R}\},
\]
\[
\operatorname*{Graph}T\cup\{(2,0.32)\}\subset\operatorname*{Graph}(\partial(\tfrac{1}{100}x^{4}))=\{(x,\tfrac{1}{25}x^{3})\mid x\in\mathbb{R}\}.
\]
However, for every $n\in\overline{3,\infty}$, we have $(0,0)\in\operatorname*{Graph}T\smallsetminus\operatorname*{Graph}(T_{n}^{++})$,
since 
\[
\operatorname*{Graph}(T_{n}^{++})\subset\operatorname*{Graph}((T_{\infty}^{+})_{n}^{+})\subset\operatorname*{Graph}((T_{\infty}^{+})_{3}^{+}),
\]
and $(0,0)\not\in\operatorname*{Graph}((T_{\infty}^{+})_{3}^{+})$,
i.e., $(0,0)$ is not $3-$m.r. to $T_{\infty}^{+}$, as evidenced
when we verify (\ref{eq:-9}) for $(x,y)=(0,0)$, $(x_{1},y_{1})=(2,0.32)$,
$(x_{2},y_{2})=(1,1)$. Indeed, in this case 
\[
(x-x_{2})y+(x_{2}-x_{1})y_{2}+(x_{1}-x)y_{1}=-0.36<0.
\]
\end{remark}

\begin{definition} \emph{Let $(X,Y,\langle\cdot,\cdot\rangle)$ be
a dual system. For $T:X\rightrightarrows Y$ and $w=(x_{0},y_{0})\in\operatorname*{Graph}T$
define }Rockafellar's antiderivative\emph{ $r_{T}^{w}:X\rightarrow\overline{\mathbb{R}}$
as } 
\begin{equation}
r_{T}^{w}(x):=\sup\{\langle x-x_{n},y_{n}\rangle+\sum_{i=0}^{n-1}\langle x_{i+1}-x_{i},y_{i}\rangle\mid n\ge1,\{(x_{i},y_{i})\}_{i\in\overline{1,n}}\subset\operatorname*{Graph}T\}.\label{RA}
\end{equation}
\end{definition}

It is worthwhile to recall Rockafellar's result, which characterizes
the convex subdifferential as a template for cyclically monotone operators.
For completeness and in consideration of the broader context, we include
a proof.

\begin{theorem} \label{RB} (\cite[Theorem 1, p. 500]{MR0193549})
Let $(X,Y,\langle\cdot,\cdot\rangle)$ be a dual system. Then $T:X\rightrightarrows Y$
is cyclically monotone iff there is $h:X\rightarrow\overline{\mathbb{R}}$
($h\in\Gamma_{\sigma(X,Y)}(X)$) such that $\operatorname*{Graph}T\subset\operatorname*{Graph}(\partial h)$.
In this case, 
\begin{equation}
\forall w=(x_{0},y_{0})\in\operatorname*{Graph}T,\ r_{T}^{w}(x_{0})=0,\ \operatorname*{Graph}T\subset\operatorname*{Graph}(\partial r_{T}^{w}),\label{eq:-23}
\end{equation}
and 
\begin{equation}
\begin{aligned}r_{T}^{w}(x) & =\min\{h(x)\mid h:X\rightarrow\overline{\mathbb{R}},\ h(x_{0})=0,\ \operatorname*{Graph}T\subset\operatorname*{Graph}(\partial h)\}\\
 & =\min\{h(x)\mid h\in\Gamma_{\sigma(X,Y)}(X),\ h(x_{0})=0,\ \operatorname*{Graph}T\subset\operatorname*{Graph}(\partial h)\}.
\end{aligned}
\label{eq:-16}
\end{equation}
Here ``$\partial$'' represents the convex subdifferential taken
with respect to $(X,Y,\langle\cdot,\cdot\rangle)$. \end{theorem}

\begin{proof} Assume that $\operatorname*{Graph}T\subset\operatorname*{Graph}(\partial h)$
for some $h:X\rightarrow\overline{\mathbb{R}}$. For every integer
$n\ge2$ and $\{(x_{i},y_{i})\}_{i\in\overline{1,n}}\subset\operatorname*{Graph}T$,
with $x_{n+1}=x_{1}$, we have 
\[
\sum_{i=1}^{n}\langle x_{i}-x_{i+1},y_{i}\rangle\ge\sum_{i=1}^{n}h(x_{i})-h(x_{i+1})=0,
\]
that is, $T$ is cyclically monotone.

Directly, let $T$ be non-empty cyclically monotone and let $w=(x_{0},y_{0})\in\operatorname*{Graph}T$.
For every $(x,y)\in\operatorname*{Graph}T$, $n\ge1$, $\{(x_{i},y_{i})\}_{i\in\overline{1,n}}\subset\operatorname*{Graph}T$
we have from the cyclic monotonicity of $T$ on the family $\{(x_{i},y_{i})\}_{i\in\overline{0,n+1}}\subset\operatorname*{Graph}T$,
where $(x_{n+1},y_{n+1})=(x,y)$ 
\[
\sum_{i=0}^{n-1}\langle x_{i}-x_{i+1},y_{i}\rangle+\langle x_{n}-x,y_{n}\rangle+\langle x-x_{0},y\rangle\ge0
\]

\vspace{-0.85cm}
 
\begin{equation}
\langle x-x_{n},y_{n}\rangle+\sum_{i=0}^{n-1}\langle x_{i+1}-x_{i},y_{i}\rangle\le\langle x-x_{0},y\rangle.\label{eq:*}
\end{equation}
After passing to supremum in (\ref{eq:*}) over $n\ge1$, $\{(x_{i},y_{i})\}_{i\in\overline{1,n}}\subset\operatorname*{Graph}T$,
we get that, for every $(x,y)\in\operatorname*{Graph}T$, $r_{T}^{w}(x)\le\langle x-x_{0},y\rangle<\infty$.
This shows that $D(T)\subset\operatorname*{dom}(r_{T}^{w})$ and $r_{T}^{w}(x_{0})\le0$.

Take $n=1$ and $(x_{1},y_{1})=(x_{0},y_{0})$ in the definition of
$r_{T}^{w}$ to find $r_{T}^{w}(x)\ge\langle x-x_{0},y_{0}\rangle$,
from which $r_{T}^{w}(x_{0})\ge0$, and so, $r_{T}^{w}(x_{0})=0$.
Note that $r_{T}^{w}\in\Gamma_{\sigma(X,Y)}(X)$ since it is a supremum
of continuous affine functions.

Arbitrarily fix $x\in X$ and $(\overline{x},\overline{y})\in\operatorname*{Graph}T$.
For every $n\ge1$, $\{(x_{i},y_{i})\}_{i\in\overline{1,n}}\subset\operatorname*{Graph}T$,
from the definition of $r_{T}^{w}$ used for $\{(x_{i},y_{i})\}_{i\in\overline{1,n+1}}\subset\operatorname*{Graph}T$,
where $(x_{n+1},y_{n+1})=(\overline{x},\overline{y})$, we find 
\[
r_{T}^{w}(x)\ge\langle x-\overline{x},\overline{y}\rangle+\langle\overline{x}-x_{n},y_{n}\rangle+\sum_{i=0}^{n-1}\langle x_{i+1}-x_{i},y_{i}\rangle.
\]

Pass to supremum over $n\ge1$, $\{(x_{i},y_{i})\}_{i\in\overline{1,n}}\subset\operatorname*{Graph}T$
to get 
\[
r_{T}^{w}(x)\ge\langle x-\overline{x},\overline{y}\rangle+r_{T}^{w}(\overline{x}),
\]
that is, $(\overline{x},\overline{y})\in\operatorname*{Graph}(\partial r_{T}^{w})$.
Therefore $\operatorname*{Graph}T\subset\operatorname*{Graph}(\partial r_{T}^{w})$.

Let $h:X\rightarrow\overline{\mathbb{R}}$ be such that $h(x_{0})=0$
and $\operatorname*{Graph}T\subset\operatorname*{Graph}(\partial h)$.

\noindent For every $\{(x_{i},y_{i})\}_{i\in\overline{1,n}}\subset\operatorname*{Graph}T$
we have $\langle x_{i+1}-x_{i},y_{i}\rangle\le h(x_{i+1})-h(x_{i})$,
$i\in\overline{1,n}$ and so 
\[
\sum_{i=0}^{n-1}\langle x_{i+1}-x_{i},y_{i}\rangle\le\sum_{i=0}^{n-1}[h(x_{i+1})-h(x_{i})]=h(x_{n}).
\]
Together with $\langle x-x_{n},y_{n}\rangle\le h(x)-h(x_{n})$ one
gets that for every $x\in X$ 
\[
\langle x-x_{n},y_{n}\rangle+\sum_{i=0}^{n-1}\langle x_{i+1}-x_{i},y_{i}\rangle\le h(x).
\]

Hence $r_{T}^{w}\le h$. \end{proof}

\begin{theorem} \label{RF} Let $(X,Y,\langle\cdot,\cdot\rangle)$
be a dual system, let $T:X\rightrightarrows Y$, and let $w=(x_{0},y_{0})\in\operatorname*{Graph}T$.
Then 
\begin{equation}
r_{T}^{(x_{0},y_{0})}(x)=\overset{{\scriptscriptstyle \langle\infty\rangle}}{\varphi}_{T}(x,y_{0})-\langle x_{0},y_{0}\rangle.\label{I}
\end{equation}
\end{theorem}

\begin{proof} Relation (\ref{I}) is an immediate consequence of
(\ref{eq:-3}), (\ref{RA}). \end{proof}

\begin{corollary} Let $(X,Y,\langle\cdot,\cdot\rangle)$ be a dual
system and let $T:X\rightrightarrows Y$ be cyclically monotone. Then,
\begin{equation}
\forall x_{0}\in D(T),\ \forall y_{0},y_{0}'\in T(x_{0}),\ r_{T}^{(x_{0},y_{0})}=r_{T}^{(x_{0},y'_{0})}.\label{eq:-38}
\end{equation}
Equivalently, 
\begin{equation}
\forall x_{0}\in D(T),\ \forall y_{0},y_{0}'\in T(x_{0}),\ \overset{{\scriptscriptstyle \langle\infty\rangle}}{\varphi}_{T}(x,y_{0})-\langle x_{0},y_{0}\rangle=\overset{{\scriptscriptstyle \langle\infty\rangle}}{\varphi}_{T}(x,y'_{0})-\langle x_{0},y'_{0}\rangle.\label{eq:-40}
\end{equation}

For every $h:X\rightarrow\overline{\mathbb{R}}$ such that $\operatorname*{Graph}T\subset\operatorname*{Graph}(\partial h)$
we have 
\begin{equation}
\forall x\in X,\ \forall(x_{0},y_{0})\in\operatorname*{Graph}T,\ h(x)\ge h(x_{0})+r_{T}^{(x_{0},y_{0})}(x).\label{eq:-39}
\end{equation}
\end{corollary}

\begin{proof} Relations (\ref{eq:-38}), (\ref{eq:-39}) are direct
consequences of (\ref{eq:-16}).

According to (\ref{I}), (\ref{eq:-40}) is equivalent to (\ref{eq:-38}).
\end{proof}

\begin{theorem} \label{E} Let $(X,Y,\langle\cdot,\cdot\rangle)$
be a dual system and let $T:X\rightrightarrows Y$. Define the operator
$K_{T}:D(K_{T})=D(T)\times R(T)\subset Z=X\times Y\rightrightarrows Z$
as $K_{T}(x,y):=T^{-1}(y)\times T(x)$, $(x,y)\in Z$, or, equivalently
\begin{equation}
((x,y),(u,v))\in\operatorname*{Graph}K_{T}\ \Leftrightarrow\ (x,v),(u,y)\in\operatorname*{Graph}T.\label{eq:-10}
\end{equation}

Then $T$ is cyclically monotone in $(X,Y,\langle\cdot,\cdot\rangle)$
iff $K_{T}$ is cyclically monotone in $(Z,Z,\cdot)$ iff $D(T)\times R(T)\subset\operatorname*{dom}\overset{{\scriptscriptstyle \langle\infty\rangle}}{\varphi}_{T}$
iff

\begin{equation}
\operatorname*{Graph}K_{T}\subset\operatorname*{Graph}(\partial\overset{{\scriptscriptstyle \langle\infty\rangle}}{\varphi}_{T})\label{eq:-22}
\end{equation}
iff 
\begin{equation}
\operatorname*{Graph}K_{T}\subset\operatorname*{Graph}(\partial\overset{{\scriptscriptstyle \langle\infty\rangle}}{\psi}_{T}).\label{eq:-31}
\end{equation}
Here ``$\partial$'' represents the convex subdifferential taken
with respect to $(Z,Z,\cdot)$. \end{theorem}

\begin{proof} We may assume that $\operatorname*{Graph}T\not=\emptyset$,
otherwise the argument is straightforward.

Fix $((x,y),(u,v))\in\operatorname*{Graph}K_{T}$, that is, fix $(x,v),(u,y)\in\operatorname*{Graph}T$.
For every $n\ge1$, $\{(x_{i},y_{i})\}_{i\in\overline{1,n}}\subset\operatorname*{Graph}T$,
use (\ref{eq:-2}), the definition of $\overset{{\scriptscriptstyle \langle n+2\rangle}}{\varphi}_{T}$,
on the pairs 
\[
(x,v),(x_{n},y_{n}),\ldots,(x_{1},y_{1}),(u,y)
\]
to get that, for every $(a,b)\in Z$, 
\[
\overset{{\scriptscriptstyle \langle n+2\rangle}}{\varphi}_{T}(a,b)\ge\langle a-x,v\rangle+\langle x-x_{n},y_{n}\rangle+\sum_{i=1}^{n-1}\langle x_{i+1}-x_{i},y_{i}\rangle+\langle x_{1}-u,y\rangle+\langle u,b\rangle.
\]
Pass to supremum over $\{(x_{i},y_{i})\}_{i\in\overline{1,n}}\subset\operatorname*{Graph}T$
to find 
\[
\forall(a,b)\in Z,\ \overset{{\scriptscriptstyle \langle n+2\rangle}}{\varphi}_{T}(a,b)\ge\langle a-x,v\rangle+\overset{{\scriptscriptstyle \langle n\rangle}}{\varphi}_{T}(x,y)+\langle u,b-y\rangle,
\]
\begin{equation}
\forall(x,v),(u,y)\in\operatorname*{Graph}T,\ \forall(a,b)\in Z,\ \overset{{\scriptscriptstyle \langle n+2\rangle}}{\varphi}_{T}(a,b)\ge(u,v)\cdot((a,b)-(x,y))+\overset{{\scriptscriptstyle \langle n\rangle}}{\varphi}_{T}(x,y).\label{eq:-24}
\end{equation}
Pass to supremum over $n\ge1$ to obtain 
\begin{equation}
\forall(x,v),(u,y)\in\operatorname*{Graph}T,\ \forall(a,b)\in Z,\ \overset{{\scriptscriptstyle \langle\infty\rangle}}{\varphi}_{T}(a,b)\ge(u,v)\cdot((a,b)-(x,y))+\overset{{\scriptscriptstyle \langle\infty\rangle}}{\varphi}_{T}(x,y).\label{eq:-25}
\end{equation}

From (\ref{eq:-25}), 
\begin{equation}
\begin{aligned}D(T)\times R(T)\subset\operatorname*{dom}\overset{{\scriptscriptstyle \langle\infty\rangle}}{\varphi}_{T}\  & \Leftrightarrow\ \forall(x,v),(u,y)\in\operatorname*{Graph}T,\ (u,v)\in\partial\overset{{\scriptscriptstyle \langle\infty\rangle}}{\varphi}_{T}(x,y)\\
 & \Leftrightarrow\ \operatorname*{Graph}K_{T}\subset\operatorname*{Graph}(\partial\overset{{\scriptscriptstyle \langle\infty\rangle}}{\varphi}_{T}).
\end{aligned}
\label{eq:-26}
\end{equation}

Conditions (\ref{eq:-22}) and (\ref{eq:-31}) are equivalent because
$(\partial\overset{{\scriptscriptstyle \langle\infty\rangle}}{\varphi}_{T})^{-1}=\partial\overset{{\scriptscriptstyle \langle\infty\rangle}}{\psi}_{T}$
and $(K_{T})^{-1}=K_{T}$.

Assume that $T$ is cyclically monotone. According to Theorem \ref{cm1},
$\overset{{\scriptscriptstyle \langle\infty\rangle}}{\varphi}_{T}\le\overset{{\scriptscriptstyle \langle\infty\rangle}}{\pisces}_{T}$.
In particular, $D(T)\times R(T)=\operatorname*{dom}\overset{{\scriptscriptstyle \langle\infty\rangle}}{\pisces}_{T}\subset\operatorname*{dom}\overset{{\scriptscriptstyle \langle\infty\rangle}}{\varphi}_{T}$.
From (\ref{eq:-26}), $\operatorname*{Graph}K_{T}\subset\operatorname*{Graph}(\partial\overset{{\scriptscriptstyle \langle\infty\rangle}}{\varphi}_{T})$.

Assume that $\operatorname*{Graph}K_{T}\subset\operatorname*{Graph}(\partial\overset{{\scriptscriptstyle \langle\infty\rangle}}{\varphi}_{T})$.
Then $K_{T}$ is cyclically monotone, since so is $\partial\overset{{\scriptscriptstyle \langle\infty\rangle}}{\varphi}_{T}$.
Also, $D(T)\times R(T)=D(K_{T})\subset\operatorname*{dom}\overset{{\scriptscriptstyle \langle\infty\rangle}}{\varphi}_{T}$.

Pick $b=y$ in (\ref{eq:-25}) to get $v\in\partial\overset{{\scriptscriptstyle \langle\infty\rangle}}{\varphi}_{T}(\cdot,y)(x)$;
whence, according to (\ref{I}), 
\begin{equation}
\operatorname*{Graph}T\subset\operatorname*{Graph}(\partial\overset{{\scriptscriptstyle \langle\infty\rangle}}{\varphi}_{T}(\cdot,y))=\operatorname*{Graph}(\partial r_{T}^{(u,y)});\label{eq:-27}
\end{equation}
in particular $T$ is cyclically monotone.

Assume that $K_{T}$ is cyclically monotone. For every finite integer
$n\ge2$ and for every $\{(x_{i},v_{i})\}_{i\in\overline{1,n}}\cup\{(u_{i},y_{i})\}_{i\in\overline{1,n}}\subset\operatorname*{Graph}T$,
the $n-$monotonicity condition of $K_{T}$ on the set $\{((x_{i},y_{i}),(u_{i},v_{i}))\}_{i\in\overline{1,n}}\subset\operatorname*{Graph}K_{T}$
translates into 
\begin{equation}
\sum_{i=1}^{n}((x_{i},y_{i})-(x_{i+1},y_{i+1}))\cdot(u_{i},v_{i})=\sum_{i=1}^{n}\langle x_{i}-x_{i+1},v_{i}\rangle+\sum_{i=1}^{n}\langle u_{i},y_{i}-y_{i+1}\rangle\ge0,\label{eq:-15}
\end{equation}
where $(x_{n+1},y_{n+1},u_{n+1},v_{n+1})=(x_{1},y_{1},u_{1},v_{1})$.
For $i\in\overline{1,n}$ we take $(u_{i},y_{i})=(u,y)$ a fixed element
of $\operatorname*{Graph}T$ to get $\sum_{i=1}^{n}\langle x_{i}-x_{i+1},v_{i}\rangle\ge0$,
i.e., $T$ is $n-$monotone. Therefore $T$ is cyclically monotone.
\end{proof}

\begin{remark} As seen in (\ref{eq:-27}), we note that (\ref{eq:-23})
is a consequence of (\ref{eq:-22}) via (\ref{eq:-25}). More precisely,
given $T$ cyclically monotone and $w=(x_{0},y_{0})\in\operatorname*{Graph}T$,
$r_{T}^{w}(x_{0})=\overset{{\scriptscriptstyle \langle\infty\rangle}}{\varphi}_{T}(x_{0},y_{0})-\langle x_{0},y_{0}\rangle=0$,
because, according to Theorem \ref{m-inf}, $\operatorname*{Graph}T\subset[\overset{{\scriptscriptstyle \langle\infty\rangle}}{\varphi}_{T}=c]$.
\end{remark}

\begin{remark} It has been observed in several places (see e.g. \cite[Theorem 3.4]{MR1009594}
or \cite[Theorem 3.3]{MR2734977}) that $T:X\rightrightarrows Y$
is monotone in $(X,Y,\langle\cdot,\cdot\rangle)$ iff $\Delta_{T}:=\{(z,z)\mid z\in\operatorname*{Graph}T\}$
is (cyclically) monotone in $(Z,Z,\cdot)$, mainly, because, for every
$n\ge1$ and $\{z_{1},\ldots,z_{n}\}\subset Z$, $z_{n+1}=z_{1}$,
we have 
\begin{align}
\sum_{i=1}^{n}c(z_{i}-z_{i+1}) & =\sum_{i=1}^{n}[c(z_{i})+c(z_{i+1})-z_{i}\cdot z_{i+1}]\nonumber \\
 & =\sum_{i=1}^{n}[2c(z_{i})-z_{i}\cdot z_{i+1}]=\sum_{i=1}^{n}(z_{i}-z_{i+1})\cdot z_{i}\label{id-cm}
\end{align}
Although $\Delta_{T}$ cannot differentiate between various types
of monotonicities, its extension $K_{T}$ can, namely, for $n\in\overline{2,\infty}$
\begin{equation}
T\ is\ n-monotone\ in\ (X,Y,\langle\cdot,\cdot\rangle)\ \Leftrightarrow\ K_{T}\ is\ n-monotone\ in\ (Z,Z,\cdot).\label{eq:-37}
\end{equation}
Indeed, the converse implication in (\ref{eq:-37}) is established
within the proof of Theorem \ref{E}.

Directly, if $T$ is $n-$monotone then so is $T^{-1}$. For every
$\{(x_{i},v_{i})\}_{i\in\overline{1,n}}\cup\{(u_{i},y_{i})\}_{i\in\overline{1,n}}\subset\operatorname*{Graph}T$,
with $(x_{n+1},y_{n+1},u_{n+1},v_{n+1})=(x_{1},y_{1},u_{1},v_{1})$,
$\sum_{i=1}^{n}\langle x_{i}-x_{i+1},v_{i}\rangle\ge0$ because $T$
is $n-$monotone, $\sum_{i=1}^{n}\langle u_{i},y_{i}-y_{i+1}\rangle\ge0$,
since $T^{-1}$ is $n-$monotone, and so, (\ref{eq:-15}) holds, that
is, $K_{T}$ is $n-$monotone. \end{remark}

Analogous to Theorem \ref{E}, the following variant of Theorem \ref{RB}
can be established.

\begin{theorem} \label{RBV} Let $(X,Y,\langle\cdot,\cdot\rangle)$
be a dual system. Then $T:X\rightrightarrows Y$ is cyclically monotone
iff one of the following conditions holds 
\begin{equation}
\forall(\exists)y_{0}\in R(T),\ D(T)\subset\operatorname*{dom}\overset{{\scriptscriptstyle \langle\infty\rangle}}{\varphi}_{T}(\cdot,y_{0}),\label{eq:-28}
\end{equation}
\begin{equation}
\forall(\exists)y_{0}\in R(T),\ \operatorname*{Graph}T\subset\operatorname*{Graph}(\partial\overset{{\scriptscriptstyle \langle\infty\rangle}}{\varphi}_{T}(\cdot,y_{0})),\label{eq:-29}
\end{equation}
\begin{equation}
\forall(\exists)w=(x_{0},y_{0})\in\operatorname*{Graph}T,\ \operatorname*{Graph}T\subset\operatorname*{Graph}(\partial r_{T}^{w}).\label{eq:-30}
\end{equation}
Here ``$\partial$'' represents the convex subdifferential taken
with respect to $(X,Y,\langle\cdot,\cdot\rangle)$. \end{theorem}

\begin{theorem} \label{KT} Let $(X,Y,\langle\cdot,\cdot\rangle)$
be a dual system, let $n\in\overline{1,\infty}$, let $T:X\rightrightarrows Y$,
and let $K_{T}:Z\rightrightarrows Z$ be defined as in (\ref{eq:-10}).
Then 
\begin{equation}
\overset{{\scriptscriptstyle \langle n\rangle}}{\varphi}_{K_{T}}(x,y,u,v)=\overset{{\scriptscriptstyle \langle n\rangle}}{\varphi}_{T}(x,v)+\overset{{\scriptscriptstyle \langle n\rangle}}{\varphi}_{T}(u,y),\ x,u\in X,\ y,v\in Y,\label{eq:-32}
\end{equation}
\begin{equation}
\overset{{\scriptscriptstyle \langle n\rangle}}{\psi}_{K_{T}}(x,y,u,v)=\overset{{\scriptscriptstyle \langle n\rangle}}{\psi}_{T}(x,v)+\overset{{\scriptscriptstyle \langle n\rangle}}{\psi}_{T}(u,y),\ x,u\in X,\ y,v\in Y,\label{eq:-33}
\end{equation}
\begin{equation}
r_{K_{T}}^{((x_{0},y_{0}),(u_{0},v_{0}))}(x,y)=r_{T}^{(x_{0},v_{0})}(x)+r_{T^{-1}}^{(y_{0},u_{0})}(y),\ (x,y)\in X\times Y,\ (x_{0},v_{0}),(u_{0},y_{0})\in\operatorname*{Graph}T,\label{eq:-34}
\end{equation}
and $T$ is (maximal) $n-$monotone (in $(X,Y,\langle\cdot,\cdot\rangle)$)
iff $K_{T}$ is (maximal) $n-$monotone (in $(Z,Z,\cdot)$).

If $T$ is cyclically monotone in $(X,Y,\langle\cdot,\cdot\rangle)$
then, for every $(x,y)\in X\times Y$, $(x_{0},v_{0})\in\operatorname*{Graph}T$,
$(u_{0},y_{0})\in\operatorname*{Graph}T$, 
\begin{equation}
r_{T^{-1}}^{(y_{0},u_{0})}(y)\le(r_{T}^{(x_{0},v_{0})})^{*}(y)-(r_{T}^{(x_{0},v_{0})})^{*}(y_{0}),\label{eq:-42}
\end{equation}
\begin{equation}
\overset{{\scriptscriptstyle \langle\infty\rangle}}{\varphi}_{T}(x,v_{0})+\overset{{\scriptscriptstyle \langle\infty\rangle}}{\varphi}_{T}(u_{0},y)-\langle x_{0},v_{0}\rangle-\langle u_{0},y_{0}\rangle\le\overset{{\scriptscriptstyle \langle\infty\rangle}}{\varphi}_{T}(x,y)-\overset{{\scriptscriptstyle \langle\infty\rangle}}{\varphi}_{T}(x_{0},y_{0}).\label{eq:-43}
\end{equation}
\end{theorem}

\begin{proof} For every finite integer $n\ge1$, for every $(x,y),(u,v)\in X\times Y$,
and for every $\{((x_{i},y_{i}),(u_{i},v_{i}))\}_{i\in\overline{1,n}}\subset\operatorname*{Graph}K_{T}$,
i.e., $\{(x_{i},v_{i})\}_{i\in\overline{1,n}}\cup\{(u_{i},y_{i})\}_{i\in\overline{1,n}}\subset\operatorname*{Graph}T$
we have 
\[
((x,y)-(x_{n},y_{n}))\cdot(u_{n},v_{n})+\sum_{i=1}^{n-1}((x_{i+1},y_{i+1})-(x_{i},y_{i}))\cdot(u_{i},v_{i})+(x_{1},y_{1})\cdot(u,v)=
\]
\[
=\langle x-x_{n},v_{n}\rangle+\sum_{i=1}^{n-1}\langle x_{i+1}-x_{i},v_{i}\rangle+\langle x_{1},v\rangle+\langle u_{n},y-y_{n}\rangle+\sum_{i=1}^{n-1}\langle u_{i},y_{i+1}-y_{i}\rangle+\langle u,y_{1}\rangle.
\]

After we pass to supremum over $\{(x_{i},v_{i})\}_{i\in\overline{1,n}}\cup\{(u_{i},y_{i})\}_{i\in\overline{1,n}}\subset\operatorname*{Graph}T$,
we find that 
\[
\overset{{\scriptscriptstyle \langle n\rangle}}{\varphi}_{K_{T}}(x,y,u,v)=\overset{{\scriptscriptstyle \langle n\rangle}}{\varphi}_{T}(x,v)+\overset{{\scriptscriptstyle \langle n\rangle}}{\varphi}_{T^{-1}}(y,u)=\overset{{\scriptscriptstyle \langle n\rangle}}{\varphi}_{T}(x,v)+\overset{{\scriptscriptstyle \langle n\rangle}}{\varphi}_{T}(u,y).
\]
Let $n\to\infty$ to complete (\ref{eq:-32}). Relation ((\ref{eq:-33}))
is obtained from (\ref{eq:-32}) by convex conjugation.

The equality in (\ref{eq:-34}) is similarly verified or it is yielded
by (\ref{eq:-32}), (\ref{I}), (\ref{eq:-41}).

If $n=1$ then $T$ is maximal $1-$monotone in $(X,Y,\langle\cdot,\cdot\rangle)$
iff $\operatorname*{Graph}T=X\times Y$ iff $\operatorname*{Graph}K_{T}=Z\times Z$
iff $K_{T}$ is maximal $1-$monotone in $(Z,Z,\cdot)$.

Assume that $T$ is maximal $n-$monotone in $(X,Y,\langle\cdot,\cdot\rangle)$,
for some $n\in\overline{2,\infty}$. According to Theorem \ref{fzv}
below, $[\overset{{\scriptscriptstyle \langle n-1\rangle}}{\varphi}_{T}\le c]\subset\operatorname*{Graph}T$
and $\overset{{\scriptscriptstyle \langle n-1\rangle}}{\varphi}_{T}\ge c$.
Similarly, since $T^{-1}$ is maximal $n-$monotone, $[\overset{{\scriptscriptstyle \langle n-1\rangle}}{\varphi}_{T^{-1}}\le c]\subset\operatorname*{Graph}(T^{-1})$
and $\overset{{\scriptscriptstyle \langle n-1\rangle}}{\varphi}_{T^{-1}}\ge c$.
After using (\ref{eq:-32}), (\ref{eq:-41}), and $\overset{{\scriptscriptstyle \langle n-1\rangle}}{\varphi}_{T}\ge c$,
$\overset{{\scriptscriptstyle \langle n-1\rangle}}{\varphi}_{T^{-1}}\ge c$,
one finds that $[\overset{{\scriptscriptstyle \langle n-1\rangle}}{\varphi}_{K_{T}}\le C]\subset\operatorname*{Graph}K_{T}$.
Here $C(z,w)=z\cdot w$, $z,w\in Z$. Therefore, again from Theorem
\ref{fzv}, $K_{T}$ is maximal $n-$monotone in $(Z,Z,\cdot)$.

Conversely, if $K_{T}$ is maximal $n-$monotone then $[\overset{{\scriptscriptstyle \langle n-1\rangle}}{\varphi}_{K_{T}}\le C]\subset\operatorname*{Graph}K_{T}$,
which yields, after considering (\ref{eq:-32}), that $[\overset{{\scriptscriptstyle \langle n-1\rangle}}{\varphi}_{T}\le c]\subset\operatorname*{Graph}T$,
i.e., $T$ is maximal $n-$monotone.

If $T$ is non-empty cyclically monotone then, for every $(x_{0},v_{0})\in\operatorname*{Graph}T$,
we have $\operatorname*{Graph}T\subset\operatorname*{Graph}(\partial r_{T}^{(x_{0},v_{0})})$.
Therefore $\operatorname*{Graph}K_{T}\subset\operatorname*{Graph}(\partial f)$,
where 
\[
f(x,y)=r_{T}^{(x_{0},v_{0})}(x)+(r_{T}^{(x_{0},v_{0})})^{*}(y),\ x\in X,\ y\in Y.
\]
Using (\ref{eq:-16}) for $K_{T}$ yields that, for every $x\in X$,
$y\in Y$, $(x_{0},v_{0}),(u_{0},y_{0})\in\operatorname*{Graph}T$,
\[
r_{K_{T}}^{((x_{0},y_{0}),(u_{0},v_{0}))}(x,y)=r_{T}^{(x_{0},v_{0})}(x)+r_{T^{-1}}^{(y_{0},u_{0})}(y)\le f(x,y)-f(x_{0},y_{0})
\]
\[
=r_{T}^{(x_{0},v_{0})}(x)+(r_{T}^{(x_{0},v_{0})})^{*}(y)-(r_{T}^{(x_{0},v_{0})})^{*}(y_{0}),
\]
from which (\ref{eq:-42}) follows.

Similarly, from $\operatorname*{Graph}K_{T}\subset\operatorname*{Graph}(\partial\overset{{\scriptscriptstyle \langle\infty\rangle}}{\varphi}_{T})$
we get 
\[
r_{K_{T}}^{((x_{0},y_{0}),(u_{0},v_{0}))}(x,y)=r_{T}^{(x_{0},v_{0})}(x)+r_{T^{-1}}^{(y_{0},u_{0})}(y)\le\overset{{\scriptscriptstyle \langle\infty\rangle}}{\varphi}_{T}(x,y)-\overset{{\scriptscriptstyle \langle\infty\rangle}}{\varphi}_{T}(x_{0},y_{0})
\]
which provides (\ref{eq:-43}) via (\ref{I}). \end{proof}

\begin{proposition} Let $(X,Y,\langle\cdot,\cdot\rangle)$ be a dual
system and let $h\in\Gamma_{\sigma(X,Y)}(X)$. Then $K_{\partial h}=\partial f$,
where $f(x,y):=h(x)+h^{*}(y)$, $x\in X$, $y\in Y$. \end{proposition}

\begin{proof} It suffices to note that $\operatorname*{Graph}(\partial h^{*})=\operatorname*{Graph}(\partial h)^{-1}$.
\end{proof}

\eject

\section{Maximal $n-$monotonicity}

\begin{theorem} \label{fzv} Let $(X,Y,\langle\cdot,\cdot\rangle)$
be a dual system and let $n\in\overline{2,\infty}$. The operator
$T:X\rightrightarrows Y$ is maximal $n-$monotone iff one of the
following conditions holds:

\medskip{}

\emph{(i)} $\operatorname*{Graph}T=(\supset)[\overset{{\scriptscriptstyle \langle n-1\rangle}}{\varphi}_{T}\le c]$;

\medskip{}

\emph{(ii)} $\overset{{\scriptscriptstyle \langle n-1\rangle}}{\varphi}_{T}\ge c$
and $\operatorname*{Graph}T=[\overset{{\scriptscriptstyle \langle n-1\rangle}}{\varphi}_{T}=c]$;

\medskip{}

\emph{(iii)} $\overset{{\scriptscriptstyle \langle n-1\rangle}}{\varphi}_{T}\ge c$
and $T$ is $n-$representable, which is defined as $T$ is $n-$monotone
and 
\begin{equation}
\exists f\in\Gamma(Z),\ c\le f\le\max\{\overset{\langle n-1\rangle}{\varphi}_{T},\overset{\langle n-1\rangle}{\psi}_{T}\},\ \operatorname*{Graph}T=[f=c];\label{eq:-18}
\end{equation}

\emph{(iv)} $T$ is $n-$representable and, for every $(x,y)\in X\times Y$,
condition $C_{n}(x,y)$ holds, where 
\begin{equation}
C_{2}(x,y):\ \exists(x_{1},y_{1})\in\operatorname*{Graph}T,\ \langle x-x_{1},y_{1}\rangle+\langle x_{1}-x,y\rangle\ge0,\label{eq:-44}
\end{equation}
\begin{equation}
\mathrm{For}\ 3\le n<\infty,\ C_{n}(x,y):\ \left\{ \begin{aligned} & \exists\{(x_{i},y_{i})\}_{i\in\overline{1,n-1}}\subset\operatorname*{Graph}T,\\
\langle x-x_{1},y_{1}\rangle & +\sum_{i=2}^{n-1}\langle x_{i-1}-x_{i},y_{i}\rangle+\langle x_{n-1}-x,y\rangle\ge0,
\end{aligned}
\right.\label{eq:-21}
\end{equation}
\begin{equation}
C_{\infty}(x,y):\ \exists m\in\mathbb{N},\ m\ge2,\ C_{m}(x,y)\ \mathrm{holds}.\label{eq:-45}
\end{equation}
\end{theorem}

\noindent\begin{proof} Assume that $T$ is maximal $n-$monotone.
Every $(x,y)\in[\overset{{\scriptscriptstyle \langle n-1\rangle}}{\varphi}_{T}\le c]$
is $n-$m.r. to $T$ which is $n-$monotone. Hence \emph{$\operatorname*{Graph}T\cup\{(x,y)\}$
}is $n-$monotone and contains $\operatorname*{Graph}T$, from which
$(x,y)\in\operatorname*{Graph}T$. We proved $[\overset{{\scriptscriptstyle \langle n-1\rangle}}{\varphi}_{T}\le c]\subset\operatorname*{Graph}T$.
Using Corollary \ref{cm2} or Theorem \ref{cm3}, we deduce the converse
inclusion from the $n-$monotonicity of $T$, thus completing the
proof of (i).

Suppose that (i) holds. Let $M$ be any $n-$monotone extension of
$T$. Every $(x,y)\in\operatorname*{Graph}M$ is $n-$m.r. to $T$,
because $M$ is $n-$monotone and contains $T$. From Theorem \ref{cm3}
(i), $(x,y)\in[\overset{{\scriptscriptstyle \langle n-1\rangle}}{\varphi}_{T}\le c]\subset\operatorname*{Graph}T$.
We proved that $\operatorname*{Graph}M\subset\operatorname*{Graph}T$.
Therefore $T$ is maximal $n-$monotone.

(ii) $\Rightarrow$ (i) is plain, since $\overset{{\scriptscriptstyle \langle n-1\rangle}}{\varphi}_{T}\ge c$
provides $[\overset{{\scriptscriptstyle \langle n-1\rangle}}{\varphi}_{T}\le c]=[\overset{{\scriptscriptstyle \langle n-1\rangle}}{\varphi}_{T}=c]$.

(i) $\Rightarrow$ (ii) According to (\ref{eq:-4}), $[\overset{{\scriptscriptstyle \langle n-1\rangle}}{\varphi}_{T}\le c]\subset\operatorname*{Graph}T\subset[\overset{{\scriptscriptstyle \langle n-1\rangle}}{\varphi}_{T}\ge c]$.
Hence $[\overset{{\scriptscriptstyle \langle n-1\rangle}}{\varphi}_{T}<c]=\emptyset$,
i.e, $\overset{{\scriptscriptstyle \langle n-1\rangle}}{\varphi}_{T}\ge c$,
and $\operatorname*{Graph}T=[\overset{{\scriptscriptstyle \langle n-1\rangle}}{\varphi}_{T}=c]$.

(ii) $\Rightarrow$ (iii) is straightforward from Theorem \ref{cm3}
(ii) with $f=\overset{{\scriptscriptstyle \langle n-1\rangle}}{\varphi}_{T}$.

(iii) $\Rightarrow$ (ii) Let $\overset{\langle n-1\rangle}{\gamma}_{T}:=\max\{\overset{\langle n-1\rangle}{\varphi}_{T},\overset{\langle n-1\rangle}{\psi}_{T}\}$.

Since $T$ is $n-$monotone we know from Theorem \ref{cm3} (ii) that
$\operatorname*{Graph}T\subset[\overset{{\scriptscriptstyle \langle n-1\rangle}}{\varphi}_{T}=c]$.

Conversely, let $z\in[\overset{{\scriptscriptstyle \langle n-1\rangle}}{\varphi}_{T}=c]$.
Then $T$ is non-empty, $\overset{{\scriptscriptstyle \langle n-1\rangle}}{\varphi}_{T}\in\Gamma(Z)$,
and, according to Lemma \ref{rep}, $z\in[\overset{{\scriptscriptstyle \langle n-1\rangle}}{\psi}_{T}=c]$;
whence $z\in[\overset{{\scriptscriptstyle \langle n-1\rangle}}{\gamma}_{T}=c]\subset[f=c]\subset\operatorname*{Graph}T$.

(iii) $\Rightarrow$ (iv) Assume that $n$ is finite. If $(x,y)\in\operatorname*{Graph}T$
then we take $(x_{i},y_{i})=(x,y)$, $i\in\overline{1,n-1}$, for
which 
\[
\langle x-x_{1},y_{1}\rangle+\langle x_{1}-x,y\rangle=\langle x-x_{1},y_{1}\rangle+\sum_{i=2}^{n-1}\langle x_{i-1}-x_{i},y_{i}\rangle+\langle x_{n-1}-x,y\rangle=0.
\]
If $(x,y)\not\in\operatorname*{Graph}T$ then, since $T$ is maximal
$n-$monotone, $\{(x,y)\}\cup\operatorname*{Graph}T$ is not $n-$monotone.
Hence, there is $\{(x_{i},y_{i})\}_{i\in\overline{1,n-1}}\subset\operatorname*{Graph}T$
such that 
\[
\langle x-x_{n-1},y\rangle+\langle x_{n-1}-x_{n-2},y_{n-1}\rangle+\ldots+\langle x_{2}-x_{1},y_{1}\rangle+\langle x_{1}-x,y_{1}\rangle<0.
\]
In both cases 
\[
\langle x-x_{1},y_{1}\rangle+\sum_{i=2}^{n-1}\langle x_{i-1}-x_{i},y_{i}\rangle+\langle x_{n-1}-x,y\rangle\ge0.
\]

When $n=\infty$ and $(x,y)\not\in\operatorname*{Graph}T$, $\{(x,y)\}\cup\operatorname*{Graph}T$
is not $\infty-$monotone. Therefore, for some $m\in\mathbb{N}$,
$m\ge2$, $\{(x,y)\}\cup\operatorname*{Graph}T$ is not $m-$monotone
and we proceed as above.

(iv) $\Rightarrow$ (iii) It suffices to prove that $\overset{{\scriptscriptstyle \langle n-1\rangle}}{\varphi}_{T}\ge c$.

For every $(x,y)\in X\times Y$ let $m\in\mathbb{N},m\ge2$ and $\{(x_{i},y_{i})\}_{i\in\overline{1,m-1}}\subset\operatorname*{Graph}T$
be such that (\ref{eq:-45}) holds, with $m=n$ when $n$ is finite.
According to $C_{n}(x,y)$, we have 
\[
\overset{{\scriptscriptstyle \langle n-1\rangle}}{\varphi}_{T}(x,y)\ge\overset{{\scriptscriptstyle \langle m-1\rangle}}{\varphi}_{T}(x,y)\ge\langle x-x_{1},y_{1}\rangle+\sum_{i=2}^{m-1}\langle x_{i-1}-x_{i},y_{i}\rangle+\langle x_{m-1},y\rangle\ge\langle x,y\rangle.
\]
\end{proof}

\begin{remark} In the specific cases where $n\in\{2,\infty\}$, 
\begin{equation}
T\ \mathit{is}\ n-\mathit{representable}\ \mathit{iff}\ \exists f\in\Gamma(Z),\ c\le f\le\overset{\langle n-1\rangle}{\psi}_{T},\ \operatorname*{Graph}T=[f=c];\label{eq:-46}
\end{equation}
because, in this case, according to Theorems \ref{m2}, \ref{m-inf},
$T$ is $n-$monotone iff $\overset{\langle n-1\rangle}{\psi}_{T}\ge\overset{\langle n-1\rangle}{\varphi}_{T}$
iff $\overset{\langle n-1\rangle}{\psi}_{T}\ge c$. In particular
\begin{equation}
\forall n\in\{2,\infty\},\ (T\ \mathit{is}\ n-\mathit{representable}\ \Rightarrow\ T\ \mathit{is}\ n-\mathit{monotone}).\label{eq:-20}
\end{equation}

Because $f\in\Gamma(Z)$, we can relax $f\le\overset{\langle n-1\rangle}{\psi_{T}}$
in (\ref{eq:-46}) to $f\le\overset{\langle n-1\rangle}{\pisces_{T}}$,
and $\operatorname*{Graph}T=[f=c]$ to $\operatorname*{Graph}T\supset[f=c]$
because from \emph{$c\le f\le\overset{\langle n-1\rangle}{\psi}_{T}\le\overset{\langle n-1\rangle}{\pisces}_{T}\le\overset{\langle1\rangle}{\pisces}_{T}=c_{T}$
}we know that $\operatorname*{Graph}T\subset[f=c]$. More precisely,
\begin{equation}
\forall n\in\{2,\infty\},\ (T\ \mathit{is}\ n-\mathit{representable}\ \mathit{iff}\ \exists f\in\Gamma(Z),\ c\le f\le\overset{\langle n-1\rangle}{\pisces}_{T},\ \operatorname*{Graph}T\supset[f=c]).\label{eq:-84}
\end{equation}

The $2-$representability definition aligns with that in \cite{MR2207807},
i.e., 
\begin{equation}
T\ \mathit{is}\ 2-\mathit{representable}\ \mathit{iff}\ \exists f\in\Gamma(Z),\ f\ge c,\ \operatorname*{Graph}T=[f=c],\label{eq:-74}
\end{equation}
since, from $\operatorname*{Graph}T\subset[f=c]$, we obtain $f\le c_{T}=\overset{\langle1\rangle}{\pisces}_{T}$.
Consequently, the condition $f\le\overset{\langle1\rangle}{\psi}_{T}$
is redundant.

Note also that, for $n\in\overline{2,\infty}$ 
\begin{equation}
T\ \mathit{is}\ \infty-\mathit{representable}\Rightarrow T\ \mathit{is}\ n-\mathit{representable}\Rightarrow(\ref{eq:-18})\Rightarrow T\ \mathit{is}\ 2-\mathit{representable}.\label{eq:-75}
\end{equation}
Indeed, the penultimate implication is plain, and the final implication
follows from (\ref{eq:-74}). For the first implication assume that
$T$ is $\infty-$representable, i.e., $\operatorname*{Graph}T=[f=c]$,
for some $f\in\Gamma(Z)$ with $c\le f\le\overset{\langle\infty\rangle}{\psi}_{T}$.
Since $T$ is $\infty-$monotone, for every $n\in\overline{2,\infty}$,
we know that $T$ is $n-$monotone and 
\[
\overset{\langle n-1\rangle}{\varphi}_{T}\le\overset{\langle\infty\rangle}{\varphi}_{T}\le\overset{\langle\infty\rangle}{\psi}_{T}\le\overset{\langle n-1\rangle}{\psi}_{T}.
\]
Therefore $f\le\overset{\langle n-1\rangle}{\psi}_{T}=\max\{\overset{\langle n-1\rangle}{\varphi}_{T},\overset{\langle n-1\rangle}{\psi}_{T}\}$
and so $T$ is $n-$representable. \end{remark}

\begin{theorem} \label{r2inf} Let $(X,Y,\langle\cdot,\cdot\rangle)$
be a dual system, let $n\in\{2,\infty\}$, and let $T:X\rightrightarrows Y$
be such that $\operatorname*{Graph}T\neq\emptyset$. Then $T$ is
$n-$representable iff $\overset{\langle n-1\rangle}{\psi}_{T}\ge c$
and $\operatorname*{Graph}T(\supset)=[\overset{\langle n-1\rangle}{\psi}_{T}=c]$
iff $T$ is $n-$monotone and $\operatorname*{Graph}T(\supset)=[\overset{\langle n-1\rangle}{\psi}_{T}=c]$.
\end{theorem}

\noindent\begin{proof} The second equivalence is plain due to Theorems
\ref{m2}, $\ref{m-inf}$. 

\noindent The converse implication in the first equivalence is plain;
take $f=\overset{\langle n-1\rangle}{\psi}_{T}$ (see also (\ref{eq:-17})).

Assume that $T$ is $n-$representable, that is, $\operatorname*{Graph}T\supset[f=c]$,
for some $f\in\Gamma(Z)$ with $c\le f\le\overset{\langle n-1\rangle}{\psi}_{T}$.
Together with (\ref{eq:-17}) we conclude that $\operatorname*{Graph}T=[\overset{\langle n-1\rangle}{\psi}_{T}=c]$,
since 
\[
[\overset{\langle n-1\rangle}{\psi}_{T}=c]\subset[f=c](\subset)=\operatorname*{Graph}T\subset[\overset{\langle n-1\rangle}{\psi}_{T}\le c]=[\overset{\langle n-1\rangle}{\psi}_{T}=c].
\]
\end{proof}

\begin{proposition} \label{sdrep} Let $(X,Y,\langle\cdot,\cdot\rangle)$
be a dual system and let $h:X\to\overline{\mathbb{R}}$. Then 
\begin{equation}
\forall(x,y)\in X\times Y,\ \overset{\langle\infty\rangle}{\pisces}_{\partial h}(x,y)\ge(h+\iota_{D(\partial h)})(x)+(h^{*}+\iota_{R(\partial h)})(y).\label{HDHI}
\end{equation}
\begin{equation}
\forall(x,y)\in X\times Y,\ \overset{\langle\infty\rangle}{\varphi}_{\partial h}(x,y)\le(h^{*}+\iota_{R(\partial h)})^{*}(x)+(h+\iota_{D(\partial h)})^{*}(y).\label{FDHI}
\end{equation}

If, in addition, $h\in\Gamma_{\sigma(X,Y)}(X)$ then $\partial h$
is cyclically representable. \end{proposition}

\noindent\begin{proof} If $\operatorname*{Graph}(\partial h)=\emptyset$
then $\overset{\langle\infty\rangle}{\pisces}_{\partial h}=+\infty$,
$\overset{\langle\infty\rangle}{\varphi}_{\partial h}=-\infty$, and
so (\ref{HDHI}), (\ref{FDHI}) are true. 

Assume that $\operatorname*{Graph}(\partial h)\not=\emptyset$. For
every $n\in\mathbb{N},$ $n\ge3$, every $(x,y)\in D(\partial h)\times R(\partial h)$,
and every 
\[
x_{1}\in(\partial h)^{-1}(y),\ y_{n}\in\partial h(x),\ \{(x_{i},y_{i})\}_{\in\overline{2,n-1}}\subset\operatorname*{Graph}(\partial h),
\]
we have $\langle x_{1},y\rangle+{\displaystyle \sum_{i=2}^{n-1}}\langle x_{i}-x_{i-1},y_{i}\rangle+\langle x-x_{n-1},y_{n}\rangle\ge$

\vspace{-.3cm}
\[
\begin{aligned}\ge h(x_{1})+h^{*}(y)+\sum_{i=2}^{n-1}[h(x_{i})-h(x_{i-1})] & +h(x)+h^{*}(y_{n})-\langle x_{n-1},y_{n}\rangle\\
=h(x_{1})+h^{*}(y)+[h(x_{n-1})-h(x_{1})] & +h(x)+h^{*}(y_{n})-\langle x_{n-1},y_{n}\rangle\\
=h(x)+h^{*}(y)+h(x_{n-1})+h^{*}(y_{n}) & -\langle x_{n-1},y_{n}\rangle\ge h(x)+h^{*}(y).
\end{aligned}
\]

Pass to infimum over $x_{1}\in(\partial h)^{-1}(y),\ y_{n}\in\partial h(x),\ \{(x_{i},y_{i})\}_{\in\overline{2,n-1}}\subset\operatorname*{Graph}(\partial h)$
to get 
\[
\forall(x,y)\in D(\partial h)\times R(\partial h),\ \forall n\in\mathbb{N},\ n\ge3,\ h(x)+h^{*}(y)\le\overset{\langle n\rangle}{\pisces}_{\partial h}(x,y),
\]
\[
\forall(x,y)\in D(\partial h)\times R(\partial h),\ h(x)+h^{*}(y)\le\inf_{n\ge3}\overset{\langle n\rangle}{\pisces}_{\partial h}(x,y)=\overset{\langle\infty\rangle}{\pisces}_{\partial h}(x,y).
\]
Since $\operatorname*{dom}(\overset{\langle\infty\rangle}{\pisces}_{\partial h})=D(\partial h)\times R(\partial h)$,
(\ref{HDHI}) holds; (\ref{FDHI}) is obtained by convex conjugation. 

If, in addition, $h\in\Gamma_{\sigma(X,Y)}(X)$, then let $f(x,y):=h(x)+h^{*}(y)$,
$x\in X$, $y\in Y$. Clearly, $f\in\Gamma(Z)$, $\operatorname*{Graph}(\partial h)=[f=c]$,
and, according to (\ref{HDHI}), $c\le f\le\overset{\langle\infty\rangle}{\pisces}_{\partial h}$.
Hence $\partial h$ is $\infty-$representable (see (\ref{eq:-84})).
\end{proof}

\begin{theorem} \label{cmax} Let $(X,Y,\langle\cdot,\cdot\rangle)$
be a dual system and let $T:X\rightrightarrows Y$ be cyclically monotone.
The following are equivalent:

\medskip{}

\noindent\emph{(i)} $T$ is maximal cyclically monotone,

\medskip{}

\noindent\emph{(ii)} $T=\partial h$, for some $h\in\Gamma_{\sigma(X,Y)}(X)$,
and, for every $(x,y)\in X\times Y$, there are $n\in\mathbb{N}$,
$n\ge2$, and $\{(x_{i},y_{i})\}_{i\in\overline{1,n-1}}\subset\operatorname*{Graph}T$
such that 
\begin{equation}
\langle x-x_{1},y_{1}\rangle+\sum_{i=2}^{n-1}\langle x_{i-1}-x_{i},y_{i}\rangle+\langle x_{n-1}-x,y\rangle\ge0,\label{eq:-296}
\end{equation}
where for $n=2$ the previous inequality reads $\langle x-x_{1},y_{1}\rangle+\langle x_{1}-x,y\rangle\ge0$.

\medskip{}

\noindent\emph{(iii)} $T=\partial h$, for some $h\in\Gamma_{\sigma(X,Y)}(X)$,
and $\overset{{\scriptscriptstyle \langle\infty\rangle}}{\varphi}_{T}\ge c$. 

\medskip

\noindent\emph{(iv)} $T=\partial h$, for some $h\in\Gamma_{\sigma(X,Y)}(X)$,
and $T$ admits a unique maximal cyclically monotone extension.\end{theorem}

\noindent\begin{proof} (i) $\Rightarrow$ (iv) is plain; (i) $\Rightarrow$
(ii) is a consequence of Theorems \ref{RB}, \ref{fzv}.

(ii) $\Rightarrow$ (iii) and (iii) $\Rightarrow$ (i) follow from
Proposition \ref{sdrep}, Theorem \ref{fzv}. 

(vi) $\Rightarrow$ (iii) It suffices to prove that $\overset{{\scriptscriptstyle \langle\infty\rangle}}{\varphi}_{T}\ge c$.
Assume by contradiction that there exists $z_{0}\in[\overset{{\scriptscriptstyle \langle\infty\rangle}}{\varphi}_{T}<c]$.
Let $f(x,y)=h(x)+h^{*}(y)$, $(x,y)\in Z$. Then $f\in\Gamma(Z)$,
$f=f^{\square}\ge c$, and  $f\ge\overset{{\scriptscriptstyle \langle\infty\rangle}}{\varphi}_{T}$,
since (\ref{HDHI}) provides $f=f^{\square}\le\overset{{\scriptscriptstyle \langle\infty\rangle}}{\pisces}_{T}$. 

Assume that $M$ is the unique maximal $\infty-$monotone extension
of $\operatorname*{Graph}T$. According to Theorem \ref{cm3} (i),
every $z\in[\overset{{\scriptscriptstyle \langle\infty\rangle}}{\varphi}_{T}\le c]$
is $\infty-$m.r. to $T$ which is $\infty-$monotone. Hence $\operatorname*{Graph}T\cup\{z\}$
is $\infty-$monotone. From Zorn's Lemma, $\operatorname*{Graph}T\cup\{z\}$
admits a maximal $\infty-$monotone extension which is also a $\infty-$monotone
extension of $\operatorname*{Graph}T$, so, every such extension coincides
with $M$. In particular $z\in M$. We proved that $[\overset{{\scriptscriptstyle \langle\infty\rangle}}{\varphi}_{T}\le c]\subset M$.
In particular $[\overset{{\scriptscriptstyle \langle\infty\rangle}}{\varphi}_{T}\le c]$
is ($2-$)monotone. 

We claim that 
\[
\operatorname*{dom}f\subset\{z_{0}\}^{+}:=\{z\in Z\mid c(z-z_{0})\geq0\}.
\]
Indeed, if there is a $z\in\operatorname*{dom}f$ $(\subset\operatorname*{dom}\overset{{\scriptscriptstyle \langle\infty\rangle}}{\varphi}_{T})$
such that $c(z-z_{0})<0$ then $z\not\in[\overset{{\scriptscriptstyle \langle\infty\rangle}}{\varphi}_{T}\le c]$,
since $[\overset{{\scriptscriptstyle \langle\infty\rangle}}{\varphi}_{T}\le c]$
is monotone and $z_{0}\in[\overset{{\scriptscriptstyle \langle\infty\rangle}}{\varphi}_{T}\le c]$.
Hence $z\in[\overset{{\scriptscriptstyle \langle\infty\rangle}}{\varphi}_{T}>c]\cap\operatorname*{dom}\overset{{\scriptscriptstyle \langle\infty\rangle}}{\varphi}_{T}$.
From the continuity of $\alpha:[0,1]\to\mathbb{R},$ $\alpha(t):=(\overset{{\scriptscriptstyle \langle\infty\rangle}}{\varphi}_{T}-c)(tz+(1-t)z_{0})$
together with $\alpha(0)\alpha(1)<0$, we infer that there exists
$s\in(0,1)\cap[\alpha=0]$, that is, $w:=sz+(1-s)z_{0}\in[\overset{{\scriptscriptstyle \langle\infty\rangle}}{\varphi}_{T}=c]$.
It follows that $c(w-z_{0})=s^{2}c(z-z_{0})<0$, which contradicts
the monotonicity of $[\overset{{\scriptscriptstyle \langle\infty\rangle}}{\varphi}_{T}\le c]$.
Our claim is proved. 

Therefore, for every $z\in\operatorname*{dom}f$ 
\[
f(z)\ge c(z)\ge c(z)-c(z-z_{0})=z\cdot z_{0}-c(z_{0})\ \Rightarrow\ c(z_{0})\ge z\cdot z_{0}-f(z).
\]
We find that $z_{0}\in[f^{\square}=c]=\operatorname*{Graph}T\subset[\overset{{\scriptscriptstyle \langle\infty\rangle}}{\varphi}_{T}=c]$
a clear contradiction. Hence $\overset{{\scriptscriptstyle \langle\infty\rangle}}{\varphi}_{T}\ge c$.
\end{proof}

\begin{remark} Note the power of the condition $T=\partial h$, for
some $h\in\Gamma_{\sigma(X,Y)}(X)$. In its presence the argument
in the proof of Theorem \ref{cmax} actually shows that 
\begin{equation}
[\overset{{\scriptscriptstyle \langle\infty\rangle}}{\varphi}_{T}\le c]\ \mathit{monotone}\ \Rightarrow\ \operatorname*{Graph}T=[\overset{{\scriptscriptstyle \langle\infty\rangle}}{\varphi}_{T}\le c]=[\overset{{\scriptscriptstyle \langle\infty\rangle}}{\varphi}_{T}=c]\ \mathit{is\ maximal}\ \infty-\mathit{monotone}.\label{eq:-85}
\end{equation}

\end{remark}

\begin{example} Let $T:D(T)=[0,+\infty)\subset\mathbb{R}\to\mathbb{R}$,
$T(x)=0$. For every $n\in\overline{2,\infty}$, $x,y\in\mathbb{R}$
\[
\overset{{\scriptscriptstyle \langle n\rangle}}{\varphi}_{T}(x,y)=\sup_{u\ge0}uy=\iota_{(-\infty,0]}(y):=\left\{ \begin{array}{cc}
0, & y\le0\\
+\infty, & y>0
\end{array}\right.,
\]
\[
\overset{{\scriptscriptstyle \langle n\rangle}}{\psi}_{T}(x,y)=(\overset{{\scriptscriptstyle \langle n\rangle}}{\varphi}_{T})^{\square}(x,y)=\sup_{v\le0,u\in\mathbb{R}}(xv+uy)=\left\{ \begin{array}{cc}
\iota_{[0,+\infty)}(x), & y=0\\
+\infty, & y\neq0
\end{array}\right.,
\]
$T$ is $\infty-$representable, since $\overset{\langle\infty\rangle}{\psi}_{T}\ge c$
and $\operatorname*{Graph}T=[\overset{\langle\infty\rangle}{\psi}_{T}=c]$,
and $T$ is not maximal $\infty-$monotone, because $\overset{\langle\infty\rangle}{\varphi}_{T}\not\ge c$.
Consequently, $T$ is neither a maximal monotone operator nor the
subdifferential of any proper, convex, and lower semicontinuous function.
\end{example}

\eject

Let $(X,Y,\langle\cdot,\cdot\rangle)$ be a dual system and let $h:X\to\overline{\mathbb{R}}$.
To $h$ we associate
\begin{equation}
\begin{aligned}h^{\cup}:X\to\overline{\mathbb{R}},\ h^{\cup}(x): & =(h^{*}+\iota_{R(\partial h)})^{*}(x)=\sup\{\langle x,y\rangle-h^{*}(y)\mid y\in R(\partial h)\}\\
 & =\sup\{\langle x-u,y\rangle+h(u)\mid(u,y)\in\operatorname*{Graph}(\partial h)\},
\end{aligned}
\label{eq:-47}
\end{equation}

\begin{equation}
\begin{aligned}h^{\circledast}:Y\to\overline{\mathbb{R}},\ h^{\circledast}(y): & =(h+\iota_{D(\partial h)})^{*}(y)=\sup\{\langle x,y\rangle-h(x)\mid x\in D(\partial h)\}\\
 & =\sup\{\langle x,y-v\rangle+h^{*}(v)\mid(x,v)\in\operatorname*{Graph}(\partial h)\},
\end{aligned}
\label{eq:-86}
\end{equation}

\begin{proposition} \label{UE} Let $(X,Y,\langle\cdot,\cdot\rangle)$
be a dual system and let $h:X\to\overline{\mathbb{R}}$. Then 

\medskip{}

\emph{(i)} $h^{\cup}\le h^{**}$, $h^{\circledast}\le h^{*}$; for
$x\in D(\partial h)$, $h^{\cup}(x)=h^{**}(x)=h(x)=h^{\circledast*}(x)$,
and, for $y\in R(\partial h)$, $h^{\cup*}(y)=h^{*}(y)=h^{\circledast}(y)$, 

\medskip{}

\emph{(ii)} $h^{\cup}$ is proper iff $\operatorname*{Graph}(\partial h)\neq\emptyset$
iff $h^{\circledast}$ is proper. In this case $h^{\cup}\in\Gamma_{\sigma(X,Y)}(X)$,
$h^{\circledast}\in\Gamma_{\sigma(Y,X)}(Y)$.

\medskip{}

\emph{(iii)} For every $n\in\mathbb{N}$, $n\ge1$, $h^{n\cup}=h^{\cup}$,
where 
\begin{equation}
\begin{aligned}h^{n\cup}(x):=\sup\{\langle x-x_{1},y_{1}\rangle+\langle x_{1}-x_{2},y_{2}\rangle+\ldots & +\langle x_{n-1}-x_{n},y_{n}\rangle+h(x_{n})\\
 & \mid\{(x_{i},y_{i})\}_{i\in\overline{1,n}}\subset\operatorname*{Graph}(\partial h)\}.
\end{aligned}
\label{nU}
\end{equation}

\emph{(iv)} $\operatorname*{Graph}(\partial h)\subset\operatorname*{Graph}(\partial h^{\cup})\cap\operatorname*{Graph}(\partial h^{\circledast})^{-1}.$ 

\medskip

Here $h^{\circledast*}:=(h^{\circledast})^{*}$, $h^{\cup*}:=(h^{\cup})^{*}$.
\end{proposition}

\begin{proof} (i) $h^{\cup}\le h^{**}$ follows directly from $h^{*}+\iota_{R(\partial h)}\ge h^{*}$
and the definition of $\partial h$. Similarly, $h+\iota_{D(\partial h)}\ge h$
implies that $h^{\circledast}\le h^{*}$. 

If $x\in D(\partial h)$ then pick $u=x$ in the definition of $h^{\cup}(x)$
to get $h^{\cup}(x)=h(x)$. Together with $h^{\cup}\le h^{**}\le h^{\circledast*}=(h+\iota_{D(\partial h)})^{**}\le h+\iota_{D(\partial h)}$
that yields $h^{\cup}(x)=h^{**}(x)=h(x)=h^{\circledast*}(x)$. In
a similar fashion, for $y\in R(\partial h)$, $h^{*}\in\Gamma(Y)$,
because $\operatorname*{Graph}(\partial h)^{-1}\subset\operatorname*{Graph}(\partial h^{*})$,
$h^{\circledast}(y)=h^{*}(y)$, after we take $v=y$ in the definition
of $h^{\circledast}$, and $h^{\circledast}\le h^{*}=h^{***}\le h^{\cup*}=(h^{*}+\iota_{R(\partial h)})^{**}\le h^{*}+\iota_{R(\partial h)}$
provides $h^{\cup*}(y)=h^{*}(y)=h^{\circledast}(y)$. 

(ii) If $\operatorname*{Graph}(\partial h)=\emptyset$ then $h^{\cup}=-\infty$,
$h^{\circledast}=-\infty$. Conversely, if $\operatorname*{Graph}(\partial h)\not=\emptyset$
then $h$ is proper, $h^{*}\in\Gamma(Y)$. Using $h^{\cup}\le h$
we see that $h^{\cup}$ is proper. Similarly, from $h^{\circledast}\le h^{*}$
we get that $h^{\circledast}$ is proper. 

(iii) Let $(u,y)\in\operatorname*{Graph}(\partial h).$ If we pick
$(x_{i},y_{i})=(u,y)$, $i\in\overline{1,n}$, in the definition of
$h^{n\cup}$ and pass to supremum over $(u,y)\in\operatorname*{Graph}(\partial h)$,
then we obtain $h^{n\cup}\ge h^{\cup}$.

Conversely, for every $\{(x_{i},y_{i})\}_{i\in\overline{1,n}}\subset\operatorname*{Graph}(\partial h)$
\[
\sum_{i=1}^{n-1}\langle x_{i}-x_{i+1},y_{i+1}\rangle\le\sum_{i=1}^{n-1}h(x_{i})-h(x_{i+1})=h(x_{1})-h(x_{n}),
\]
from which $h^{n\cup}\le h^{\cup}$.

(iv) For every $(x,y)\in\operatorname*{Graph}(\partial h)$, $w\in X$,
we have $h^{\cup}(x)=h(x)$ and, from the definition of $h^{\cup}(w)$,
$\langle w-x,y\rangle+h(x)\le h^{\cup}(w)$. Hence $(x,y)\in\operatorname*{Graph}(\partial h^{\cup})$. 

Analogously, for every $(x,y)\in\operatorname*{Graph}(\partial h)$,
$\zeta\in Y$, we have $h^{\circledast}(\zeta)=h^{*}(\zeta)$ and,
from the definition of $h^{\circledast}(\zeta)$, $\langle x,\zeta-y\rangle+h^{*}(y)\le h^{\circledast}(\zeta)$.
Hence $(y,x)\in\operatorname*{Graph}(\partial h^{\circledast})$.
\end{proof}

\begin{theorem} \label{G} Let $(X,Y,\langle\cdot,\cdot\rangle)$
be a dual system and let $h:X\to\overline{\mathbb{R}}$. For every
$(x,y)\in X\times Y$, $(x_{0},y_{0})\in\operatorname*{Graph}(\partial h)$
\begin{equation}
\overset{{\scriptscriptstyle \langle\infty\rangle}}{\varphi}_{\partial h}(x,y)=h^{\cup}(x)+h^{\circledast}(y),\ \overset{{\scriptscriptstyle \langle\infty\rangle}}{\psi}_{\partial h}(x,y)=h^{\circledast*}(x)+h^{\cup*}(y)\label{eq:-52}
\end{equation}
\begin{equation}
r_{\partial h}^{(x_{0},y_{0})}(x)=h^{\cup}(x)-h(x_{0}),\ r_{(\partial h)^{-1}}^{(y_{0},x_{0})}(x)=h^{\circledast}(y)-h^{*}(y_{0}).\label{eq:-67}
\end{equation}
\end{theorem}

\begin{proof} If $\operatorname*{Graph}(\partial h)=\emptyset$ then
(\ref{eq:-52}), (\ref{eq:-67}) are trivially verified with $\overset{{\scriptscriptstyle \langle\infty\rangle}}{\varphi}_{\partial h}=h^{\cup}=h^{\circledast}=r_{\partial h}^{(x_{0},y_{0})}=-\infty$,
$\overset{{\scriptscriptstyle \langle\infty\rangle}}{\psi}_{\partial h}=h^{\circledast*}=h^{\cup*}=+\infty$. 

Assume that $\operatorname*{Graph}(\partial h)\not=\emptyset$. According
to Proposition \ref{UE} (ii), $h^{\cup}\in\Gamma(X)$, $h^{\circledast}\in\Gamma(Y)$. 

Let us prove that, for every $x\in X$, $y\in Y$, 
\begin{equation}
\sup_{v_{0}\in R(\partial h)}\overset{{\scriptscriptstyle \langle\infty\rangle}}{\varphi}_{\partial h}(x,v_{0})-h^{*}(v_{0})=h^{\cup}(x),\ \sup_{u_{0}\in D(\partial h)}\overset{{\scriptscriptstyle \langle\infty\rangle}}{\varphi}_{\partial h}(u_{0},y)-h(u_{0})=h^{\circledast}(y).\label{eq:-53}
\end{equation}

According to (\ref{eq:-4}), for every $x\in X$, $v_{0}\in R(\partial h)$,
$\overset{{\scriptscriptstyle \langle\infty\rangle}}{\varphi}_{\partial h}(x,v_{0})\ge\langle x,v_{0}\rangle$,
and so $\sup_{v_{0}\in R(\partial h)}\overset{{\scriptscriptstyle \langle\infty\rangle}}{\varphi}_{\partial h}(x,v_{0})-h^{*}(v_{0})\ge\sup_{v_{0}\in R(\partial h)}\langle x,v_{0}\rangle-h^{*}(v_{0})=h^{\cup}(x)$.
The converse inequality is a consequence of (\ref{FDHI}) and $h^{*}=h^{\circledast}$
on $R(\partial h)$, or $h^{\circledast}\le h^{*}$. We proceed similarly
for the second part in (\ref{eq:-53}). 

We use (\ref{eq:-43}) for $x\in X$, $(x_{0},v_{0}),(u_{0},y_{0})\in\operatorname*{Graph}(\partial h)\subset[\overset{{\scriptscriptstyle \langle\infty\rangle}}{\varphi}_{\partial h}=c]$,
and $y=y_{0}$ to get that, for every $x\in X$, $y_{0}\in R(\partial h)$,
$(x_{0},v_{0})\in\operatorname*{Graph}(\partial h)$ 
\[
\begin{aligned}\overset{{\scriptscriptstyle \langle\infty\rangle}}{\varphi}_{\partial h}(x,y_{0}) & \ge\overset{{\scriptscriptstyle \langle\infty\rangle}}{\varphi}_{\partial h}(x_{0},y_{0})+\overset{{\scriptscriptstyle \langle\infty\rangle}}{\varphi}_{\partial h}(x,v_{0})-\langle x_{0},v_{0}\rangle\\
 & =[\overset{{\scriptscriptstyle \langle\infty\rangle}}{\varphi}_{\partial h}(x,v_{0})-h^{*}(v_{0})]+[\overset{{\scriptscriptstyle \langle\infty\rangle}}{\varphi}_{\partial h}(x_{0},y_{0})-h(x_{0})].
\end{aligned}
\]
Pass to supremum over $(x_{0},v_{0})\in\operatorname*{Graph}(\partial h)$
to obtain, according to (\ref{eq:-53}), (\ref{FDHI}), that, for
every $x\in X$, $y_{0}\in R(\partial h)$, 
\begin{equation}
\overset{{\scriptscriptstyle \langle\infty\rangle}}{\varphi}_{\partial h}(x,y_{0})=h^{\cup}(x)+h^{\circledast}(y_{0})=h^{\cup}(x)+h^{*}(y_{0}).\label{eq:-55}
\end{equation}
Similarly, for every $y\in Y$, $x_{0}\in D(\partial h)$, 
\begin{equation}
\overset{{\scriptscriptstyle \langle\infty\rangle}}{\varphi}_{\partial h}(x_{0},y)=h(x_{0})+h^{\circledast}(y).\label{eq:-56}
\end{equation}

For every $(x,y)\in X\times Y$, we use (\ref{eq:-43}) for $(x_{0},v_{0})=(u_{0},y_{0})\in\operatorname*{Graph}(\partial h)$
to get 
\[
\overset{{\scriptscriptstyle \langle\infty\rangle}}{\varphi}_{\partial h}(x,y)-\langle x_{0},y_{0}\rangle=\overset{{\scriptscriptstyle \langle\infty\rangle}}{\varphi}_{\partial h}(x,y)-\overset{{\scriptscriptstyle \langle\infty\rangle}}{\varphi}_{\partial h}(x_{0},y_{0})\ge\overset{{\scriptscriptstyle \langle\infty\rangle}}{\varphi}_{\partial h}(x,v_{0})+\overset{{\scriptscriptstyle \langle\infty\rangle}}{\varphi}_{\partial h}(u_{0},y)-\langle x_{0},v_{0}\rangle-\langle u_{0},y_{0}\rangle
\]
\[
=h^{\cup}(x)+h^{*}(v_{0})+h(u_{0})+h^{\circledast}(y)-\langle x_{0},v_{0}\rangle-\langle u_{0},y_{0}\rangle=h^{\cup}(x)+h^{\circledast}(y)-h(x_{0})-h^{*}(y_{0}),
\]
so 
\[
\overset{{\scriptscriptstyle \langle\infty\rangle}}{\varphi}_{\partial h}(x,y)\ge h^{\cup}(x)+h^{\circledast}(y).
\]

The first equality in (\ref{eq:-52}) follows from (\ref{FDHI}).
The second equality in (\ref{eq:-52}) is obtained by convex conjugation
from the first part. For (\ref{eq:-67}) we use (\ref{I}), (\ref{eq:-52}),
(\ref{eq:-41}), and $h^{\cup}=h$ on $D(\partial h)$, $h^{\circledast}=h^{*}$
on $R(\partial h)$. \end{proof}

\begin{corollary} \label{RTE} Let $(X,Y,\langle\cdot,\cdot\rangle)$
be a dual system and let $h:X\to\overline{\mathbb{R}}$. Then 
\begin{equation}
\forall x_{0}\in D(\partial h),\ h^{\cup}(x)=\min\{\alpha(x)\mid\alpha:X\to\overline{\mathbb{R}},\ \alpha(x_{0})=h(x_{0}),\ \operatorname*{Graph}(\partial\alpha)\subset\operatorname*{Graph}(\partial h)\},\label{eq:-58}
\end{equation}
\begin{equation}
\forall y_{0}\in R(\partial h),\ h^{\circledast}(y)=\min\{\beta(y)\mid\beta:Y\to\overline{\mathbb{R}},\ \beta(y_{0})=h^{*}(y_{0}),\ \operatorname*{Graph}(\partial\beta)\subset\operatorname*{Graph}(\partial h)^{-1}\}.\label{eq:-48}
\end{equation}
\end{corollary}

\begin{proof} Relations (\ref{eq:-58}), (\ref{eq:-48}) are translations
of (\ref{eq:-16}) with the aid of (\ref{eq:-67}). \end{proof}

\begin{corollary} \label{ES} Let $(X,Y,\langle\cdot,\cdot\rangle)$
be a dual system and let $h,g\in\Gamma_{\sigma(X,Y)}(X)$. Then $\partial h=\partial g$
iff $\overset{{\scriptscriptstyle \langle\infty\rangle}}{\varphi}_{\partial h}=\overset{{\scriptscriptstyle \langle\infty\rangle}}{\varphi}_{\partial g}$
iff $h^{\cup}-g^{\cup}=h^{\circledast*}-g^{\circledast*}$ is a constant
function. \end{corollary}

\begin{proof} The conclusion is straightforward if $\operatorname*{Graph}(\partial h)=\emptyset$
or $\operatorname*{Graph}(\partial g)=\emptyset$; in which case $\overset{{\scriptscriptstyle \langle\infty\rangle}}{\varphi}_{\partial h}=\overset{{\scriptscriptstyle \langle\infty\rangle}}{\varphi}_{\partial g}=-\infty$
and $h^{\cup}-g^{\cup}=h^{\circledast*}-g^{\circledast*}=+\infty$. 

Assume that $\operatorname*{Graph}(\partial h)\not=\emptyset$ and
$\operatorname*{Graph}(\partial g)\not=\emptyset$; whence $h^{\cup},g^{\cup}\in\Gamma(X)$,
$h^{\circledast},g^{\circledast}\in\Gamma(Y)$. 

If $\partial h=\partial g$ then $\overset{{\scriptscriptstyle \langle\infty\rangle}}{\varphi}_{\partial h}=\overset{{\scriptscriptstyle \langle\infty\rangle}}{\varphi}_{\partial g}$.
From (\ref{eq:-52}) we get that $h^{\cup}=g^{\cup}+k$, where $k=g^{\circledast}(y_{0})-h^{\circledast}(y_{0})=g^{*}(y_{0})-h^{*}(y_{0})\in\mathbb{R}$,
for any fixed $y_{0}\in R(\partial h)=R(\partial g)$. Similarly,
$h^{\circledast}=g^{\circledast}+h$, for some $h\in\mathbb{R}$.
Again $\overset{{\scriptscriptstyle \langle\infty\rangle}}{\varphi}_{\partial h}=\overset{{\scriptscriptstyle \langle\infty\rangle}}{\varphi}_{\partial g}$
shows that $h=-k$, and so, $h^{\circledast*}-g^{\circledast*}=k$. 

Conversely, if $h^{\cup}-g^{\cup}=h^{\circledast*}-g^{\circledast*}=k\in\mathbb{R}$
then $h^{\circledast}=g^{\circledast}-k$, $\overset{{\scriptscriptstyle \langle\infty\rangle}}{\varphi}_{\partial h}=\overset{{\scriptscriptstyle \langle\infty\rangle}}{\varphi}_{\partial g}$,
followed by $\overset{{\scriptscriptstyle \langle\infty\rangle}}{\psi}_{\partial h}=\overset{{\scriptscriptstyle \langle\infty\rangle}}{\psi}_{\partial g}$.
According to Proposition \ref{sdrep}, $\partial h$, $\partial g$
are $\infty-$representable. From Theorem \ref{r2inf}, we find that
\[
\operatorname*{Graph}(\partial h)=[\overset{{\scriptscriptstyle \langle\infty\rangle}}{\psi}_{\partial h}=c]=[\overset{{\scriptscriptstyle \langle\infty\rangle}}{\psi}_{\partial g}=c]=\operatorname*{Graph}(\partial g).
\]
\end{proof}

\begin{remark} \emph{\label{Ech} Whenever }$h:X\to\overline{\mathbb{R}}$
has\emph{ }$\operatorname*{Graph}(\partial h)\neq\emptyset$ we have
\begin{equation}
h^{**}=h^{\cup},\ h^{*}=h^{\circledast}\ \Leftrightarrow\ h^{\cup}=h^{\circledast*}\ \Leftrightarrow\ h^{\cup*}=h^{\circledast}.\label{eq:-76}
\end{equation}
\emph{Indeed, for the first equivalence, if $h^{**}=h^{\cup}$ and
$h^{*}=h^{\circledast}$ then, by conjugation, $h^{\cup}=h^{\circledast*}$.
Conversely, if $h^{\cup}=h^{\circledast*}$, from $h^{\cup}\le h^{**}\le h^{\circledast*}$
we know that $h^{**}=h^{\cup}=h^{\circledast*}$ and $h^{*}=h^{\circledast}$,
since $h^{*},h^{\circledast}\in\Gamma(Y)$. The second equivalence
follows again from the biconjugate formula applied for $h^{\cup}\in\Gamma(X)$
and $h^{\circledast}\in\Gamma(Y)$. }\end{remark}

\begin{theorem} \label{F} Let $(X,Y,\langle\cdot,\cdot\rangle)$
be a dual system and let $h\in\Gamma_{\sigma(X,Y)}(X)$. The following
are equivalent:

\medskip{}

\emph{(i)} $\partial h$ is maximal cyclically monotone,

\medskip{}

\emph{(ii)} $h^{\cup}=h^{\circledast*}$, 

\medskip{}

\emph{(iii)} $\forall(x,y)\in X\times Y$, $\overset{{\scriptscriptstyle \langle\infty\rangle}}{\varphi}_{\partial h}(x,y)=h(x)+h^{*}(y)$,

\medskip{}

\emph{(iv)} $\forall(x,y)\in X\times Y$, $\forall(x_{0},v_{0})\in\operatorname*{Graph}(\partial h)$,
\[
r_{\partial h}^{(x_{0},y_{0})}(x)=h(x)-h(x_{0}),\ r_{\partial h^{*}}^{(y_{0},x_{0})}(y)=h^{*}(y)-h^{*}(y_{0}),
\]

\emph{(v)} $\forall(x,y)\in X\times Y$, $\forall x_{0}\in D(\partial h)$,
$\forall y_{0}\in R(\partial h)$, 
\[
\overset{{\scriptscriptstyle \langle\infty\rangle}}{\varphi}_{\partial h}(x,y_{0})=h(x)+h^{*}(y_{0}),\ \overset{{\scriptscriptstyle \langle\infty\rangle}}{\varphi}_{\partial h}(x_{0},y)=h(x_{0})+h^{*}(y),
\]

\emph{(vi)} $\forall(x,y)\in X\times Y$, $\exists(x_{0},y_{0})\in\operatorname*{Graph}(\partial h)$
\[
\overset{{\scriptscriptstyle \langle\infty\rangle}}{\varphi}_{\partial h}(x,y_{0})=h(x)+h^{*}(y_{0}),\ \overset{{\scriptscriptstyle \langle\infty\rangle}}{\varphi}_{\partial h}(x_{0},y)=h(x_{0})+h^{*}(y).
\]
\end{theorem}

\begin{proof} (i) $\Rightarrow$ (ii) Since $\partial h$ is maximal
$\infty-$monotone, we know that $\operatorname*{Graph}(\partial h)\neq\emptyset$,
and so, $h^{\cup}\in\Gamma(X)$, $h^{\circledast}\in\Gamma(Y)$.

For every $x\in X$ and $y=0$ we use Theorem \ref{cmax} to get $n\in\mathbb{N}$,
$n\ge2$ and $\{(x_{i},y_{i})\}_{i\in\overline{1,n-1}}\subset\operatorname*{Graph}(\partial h)$
such that 
\[
\langle x-x_{1},y_{1}\rangle+\sum_{i=2}^{n-1}\langle x_{i-1}-x_{i},y_{i}\rangle\ge0,
\]
where for $n=2$ the previous inequality reads $\langle x-x_{1},y_{1}\rangle\ge0$. 

Therefore, for every $x\in X$ 
\[
\begin{aligned}h(x)\ge h^{\cup}(x)=h^{(n-1)\cup}(x) & \ge\langle x-x_{1},y_{1}\rangle+\sum_{i=2}^{n-1}\langle x_{i-1}-x_{i},y_{i}\rangle+h(x_{n-1}\rangle\\
 & \ge h(x_{n-1})\ge\inf_{D(\partial h)}h\ge\inf_{X}h,
\end{aligned}
\]
and so, $\inf_{X}h=\inf_{D(\partial h)}h=\inf_{X}h^{\cup}$ or $h^{\circledast}(0)=h^{\cup*}(0)$. 

But $\partial h$ maximal $\infty-$monotone implies (is equivalent
to) that $\partial(h-y)$ is maximal $\infty-$monotone, for every
(some) $y\in Y$; where $(h-y)(x):=h(x)-\langle x,y\rangle$, $x\in X$,
because $\operatorname*{Graph}(\partial(h-y))=\operatorname*{Graph}(\partial h)-(0,y)$.
Hence, for every $y\in Y$, $(h-y)^{\cup*}(0)=(h-y)^{\circledast}(0)$. 

Note that, for $x\in X$, $y,v\in Y$, $(h-y)^{\cup}(x)=h^{\cup}(x)-\langle x,y\rangle$,
$(h-y)^{\cup*}(v)=h^{\cup*}(y+v)$, $(h-y)^{\circledast}(v)=h^{\circledast}(y+v)$.
Set $v=0$ to get $h^{\cup*}=h^{\circledast}$ and, by conjugation,
$h^{\cup}=h^{\circledast*}$. 

For an alternative argument for this implication, we use Theorem \ref{cmax}
and (\ref{eq:-81}) below. 

(ii) $\Rightarrow$ (iii) follows from (\ref{eq:-52}) (see also Remark
\ref{Ech}).

(iii) $\Rightarrow$ (i) is a direct consequence of Proposition \ref{sdrep},
and Theorem \ref{fzv} or Theorem \ref{cmax}. 

(iv) and (v) are equivalent due to (\ref{I}), $\partial h^{*}=(\partial h)^{-1}$,
and (\ref{eq:-41}). 

(iii) $\Rightarrow$ (v) $\Rightarrow$ (vi) are plain.

(vi) $\Rightarrow$ (iii) For every $(x,y)\in X\times Y$, let $(x_{0},y_{0})\in\operatorname*{Graph}(\partial h)$
be such that $\overset{{\scriptscriptstyle \langle\infty\rangle}}{\varphi}_{\partial h}(x,y_{0})=h(x)+h^{*}(y_{0})$,
$\overset{{\scriptscriptstyle \langle\infty\rangle}}{\varphi}_{\partial h}(x_{0},y)=h(x_{0})+h^{*}(y)$. 

We use (\ref{eq:-43}) for $(x_{0},v_{0})=(u_{0},y_{0})\in\operatorname*{Graph}(\partial h)$
to get 
\[
\overset{{\scriptscriptstyle \langle\infty\rangle}}{\varphi}_{\partial h}(x,y)-\langle x_{0},y_{0}\rangle=\overset{{\scriptscriptstyle \langle\infty\rangle}}{\varphi}_{\partial h}(x,y)-\overset{{\scriptscriptstyle \langle\infty\rangle}}{\varphi}_{\partial h}(x_{0},y_{0})\ge\overset{{\scriptscriptstyle \langle\infty\rangle}}{\varphi}_{\partial h}(x,v_{0})+\overset{{\scriptscriptstyle \langle\infty\rangle}}{\varphi}_{\partial h}(u_{0},y)-\langle x_{0},v_{0}\rangle-\langle u_{0},y_{0}\rangle
\]
\[
=h(x)+h^{*}(v_{0})+h(u_{0})+h^{*}(y)-\langle x_{0},v_{0}\rangle-\langle u_{0},y_{0}\rangle=h(x)+h^{*}(y)-h(x_{0})-h^{*}(y_{0}),
\]
so 
\[
\overset{{\scriptscriptstyle \langle\infty\rangle}}{\varphi}_{\partial h}(x,y)\ge h(x)+h^{*}(y).
\]
Together with (\ref{eq:-52}) and $h^{\cup}\le h$, $h^{\circledast}\le h^{*}$
we find $\overset{{\scriptscriptstyle \langle\infty\rangle}}{\varphi}_{\partial h}(x,y)=h(x)+h^{*}(y)$
and (ii). \end{proof}

\begin{corollary} Let $(X,Y,\langle\cdot,\cdot\rangle)$ be a dual
system and let $h,g\in\Gamma_{\sigma(X,Y)}(X)$ be such that $\partial h$
is maximal cyclically monotone. Then $\partial h(\subset)=\partial g$
iff $h-g$ is a constant function. \end{corollary}

\begin{proof} The converse is plain. For the direct implication,
given that $\partial h=\partial g$ are both maximal cyclically monotone,
from Theorem \ref{F} and Proposition \ref{UE} (i) or Remark \ref{Ech},
we infer that $h^{\cup}=h^{\circledast*}=h$, $g^{\cup}=g^{\circledast*}=g$,
and use Corollary \ref{ES} to conclude. \end{proof}

\strut

We end this section with some conclusive remarks concerning the maximal
cyclical monotonicity of $\partial h$, where $h:X\to\overline{\mathbb{R}}$
is a general function. 

\medskip

If $\partial h$ is maximal cyclically monotone then, since $\operatorname*{Graph}(\partial h)\subset\operatorname*{Graph}(\partial h^{**})$,
we infer that $h^{**}\in\Gamma_{\sigma(X,Y)}(X)$, $h^{*}\in\Gamma_{\sigma(Y,X)}(Y)$,
and $\partial h=\partial h^{**}$ is maximal cyclically monotone.
According to Theorem \ref{F}, $h^{**\cup}=h^{**\circledast*}=h^{**}$.
The converse is not true, since one can find $h:X\to\overline{\mathbb{R}}$
such that $\operatorname*{Graph}(\partial h)\subsetneq\operatorname*{Graph}(\partial h^{**})$
and $\partial h^{**}$ is maximal (cyclically) monotone, e.g., when
$X$ is a Banach space and $h=\iota_{C}$, where $C\subsetneq X$
is open convex. 

Note that, in all of the arguments above, the condition $h\in\Gamma_{\sigma(X,Y)}(X)$
is essential only in establishing the cyclical representability of
$\partial h$. 

Similar computations show that, in general, for $h,g:X\to\overline{\mathbb{R}}$
\begin{equation}
\overset{{\scriptscriptstyle \langle\infty\rangle}}{\varphi}_{\partial h}\ge c\ \Leftrightarrow\ h^{\cup}(\ge)=h^{\circledast*}(=h^{**})\ \Rightarrow\ \partial h^{**}\ \mathit{is\ maximal}\ \infty-\mathit{monotone}.\label{eq:-81}
\end{equation}
Indeed, from (\ref{eq:-52}), we know that $\overset{{\scriptscriptstyle \langle\infty\rangle}}{\varphi}_{\partial h}(x,y)=h^{\cup}(x)+h^{\circledast}(y)$,
$x\in X$, $y\in Y$. Therefore $\overset{{\scriptscriptstyle \langle\infty\rangle}}{\varphi}_{\partial h}\ge c$
is equivalent to $h^{\cup}\ge h^{\circledast*}$. However, $h^{\cup}\le h^{**}\le h^{\circledast*}$. 

If $\overset{{\scriptscriptstyle \langle\infty\rangle}}{\varphi}_{\partial h}\ge c$,
then by Theorem \ref{fzv}, any cyclically representable extension
of $\partial h$, such as $\partial h^{**}$, is maximal cyclically
monotone. 

Similarly, 
\begin{equation}
\begin{aligned}\overset{{\scriptscriptstyle \langle\infty\rangle}}{\varphi}_{\partial h}=\overset{{\scriptscriptstyle \langle\infty\rangle}}{\varphi}_{\partial g} & \ \Leftrightarrow\\
\Leftrightarrow\ (h^{*}+\iota_{R(\partial h)})^{*}-(g^{*}+\iota_{R(\partial g)})^{*} & =(h+\iota_{D(\partial h)})^{**}-(g+\iota_{D(\partial g)})^{**}\ \mathit{is\ constant}.
\end{aligned}
\label{eq:-80}
\end{equation}

\eject

\section{Examples }

\begin{example} (\textsc{Monotone Subset of} $\mathbb{R}^{2}$) If
$T:\mathbb{R}\rightrightarrows\mathbb{R}$ is monotone, then $T$
is cyclically monotone. \end{example}

\begin{proof} We prove by induction over $n\in\mathbb{N}$, $n\ge2$
that $T$ is $n-$monotone. Note that the base case is provided. Assume
that $T$ is $n-$monotone, for some $n\in\mathbb{N}$, $n\ge2$.

Let $\{(x_{i},y_{i})\}_{i\in\overline{1,n+1}}\subset\operatorname*{Graph}T$
be arbitrarily chosen with $x_{n+2}=x_{1}$. Pick $k\in\overline{1,n+1}$
such that $x_{k}=\min_{i\in\overline{1,n+1}}x_{i}$ and define $(x_{j}',y_{j}')=(x_{(j+k-1)\mathrm{mod}\,n},y_{(j+k-1)\mathrm{mod}\,n})$,
$j\in\overline{1,n+1}$. Then $\{(x'_{j},y'_{j})\}_{i\in\overline{1,n+1}}\subset\operatorname*{Graph}T$,
$x_{1}'=\min_{j\in\overline{1,n+1}}x'_{j}$, and, because $T$ is
monotone, $y_{1}'=\min_{j\in\overline{1,n+1}}y'_{j}$. We have $\sum_{j=2}^{n}(x_{j}'-x_{j+1}')y_{j}'+(x_{n+1}'-x_{2}')y'_{n+1}\ge0$,
since $T$ is $n-$monotone, and so 
\[
\begin{aligned}\sum_{i=1}^{n+1}(x_{i}-x_{i+1})y_{i} & =\sum_{j=1}^{n+1}(x_{j}'-x_{j+1}')y_{j}'\\
 & =(x_{1}'-x_{2}')(y_{1}'-y'_{n+1})+\sum_{j=2}^{n}(x_{j}'-x_{j+1}')y_{j}'+(x_{n+1}'-x_{2}')y'_{n+1}\ge0,
\end{aligned}
\]
because $x_{1}'\le x_{2}'$, $y_{1}'\le y'_{n+1}$. Therefore $T$
is $(n+1)-$monotone and the proof is complete. \end{proof}

\begin{theorem} \label{R} Let $(X,Y,\langle\cdot,\cdot\rangle)$
be a dual system, let $T:X\rightrightarrows Y$, and let $n\in\overline{1,\infty}$.
Then 
\begin{equation}
\begin{aligned}\overset{{\scriptscriptstyle \langle n+1\rangle}}{\pisces}_{T}(x,y) & =\negthickspace\negthickspace\negthickspace\inf_{u\in T^{-1}(y),v\in R(T)}\left\{ \overset{{\scriptscriptstyle \langle n\rangle}}{\pisces}_{T}(x,v)+\langle u,y-v\rangle\right\} \\
 & =\negthickspace\negthickspace\negthickspace\inf_{u\in D(T),v\in T(x)}\left\{ \overset{{\scriptscriptstyle \langle n\rangle}}{\pisces}_{T}(u,y)+\langle x-u,v\rangle\right\} ;
\end{aligned}
\label{eq:-61}
\end{equation}
\begin{equation}
\begin{aligned}\overset{{\scriptscriptstyle \langle n+1\rangle}}{\varphi}_{T}(x,y) & =\negthickspace\negthickspace\negthickspace\sup_{(u,v)\in\operatorname*{Graph}T}\left\{ \overset{{\scriptscriptstyle \langle n\rangle}}{\varphi}_{T}(x,v)+\langle u,y-v\rangle\right\} \\
 & =\negthickspace\negthickspace\negthickspace\sup_{(u,v)\in\operatorname*{Graph}T}\left\{ \overset{{\scriptscriptstyle \langle n\rangle}}{\varphi}_{T}(u,y)+\langle x-u,v\rangle\right\} .
\end{aligned}
\label{eq:-60}
\end{equation}
\end{theorem}

\begin{proof} For $n=1$, (\ref{eq:-61}) is directly verified. For
$n=2$ 
\[
\begin{aligned}\overset{{\scriptscriptstyle \langle3\rangle}}{\pisces}_{T}(x,y) & =\inf_{x_{1}\in T^{-1}(y),\ y_{3}\in T(x),\ (x_{2},y_{2})\in\operatorname*{Graph}T}\{\langle x_{1},y\rangle+\langle x_{2}-x_{1},y_{2}\rangle+\langle x-x_{2},y_{3}\rangle\}\\
 & =\inf_{x_{2}\in D(T),\ y_{3}\in T(x)}\left(\inf_{x_{1}\in T^{-1}(y),\ y_{2}\in T(x_{2})}\{\langle x_{1},y\rangle+\langle x_{2}-x_{1},y_{2}\rangle\}+\langle x-x_{2},y_{3}\rangle\right)\\
 & =\inf_{x_{2}\in D(T),\ y_{3}\in T(x)}\left\{ \overset{{\scriptscriptstyle \langle2\rangle}}{\pisces}_{T}(x_{2},y)+\langle x-x_{2},y_{3}\rangle\right\} .
\end{aligned}
\]
For $n\in\mathbb{N}$, $n\ge3$ the second equality in (\ref{eq:-61})
holds analogously, as follows 
\[
\begin{aligned}\overset{{\scriptscriptstyle \langle n+1\rangle}}{\pisces}_{T}(x,y)=\inf\left\{ \langle x_{1},y\rangle+\right. & \sum_{i=2}^{n}\langle x_{i}-x_{i-1},y_{i}\rangle+\langle x-x_{n},y_{n+1}\rangle\mid\\
 & \left.x_{1}\in T^{-1}(y),\ y_{n+1}\in T(x),\ (x_{i},y_{i})\}_{i\in\overline{2,n}}\subset\operatorname*{Graph}T\right\} \\
=\inf_{x_{n}\in D(T),y_{n+1}\in T(x)}\left[\inf\right.\left\{ \langle x_{1},y\rangle+\right. & \sum_{i=2}^{n-1}\langle x_{i}-x_{i-1},y_{i}\rangle+\langle x_{n}-x_{n-1},y_{n}\rangle\mid\\
 & x_{1}\in T^{-1}(y),\ y_{n}\in T(x_{n}),\ (x_{i},y_{i})\}_{i\in\overline{2,n-1}}\subset\operatorname*{Graph}\left.T\right\} \\
\left.+\langle x-x_{n},y_{n+1}\rangle\right] & =\negthickspace\negthickspace\negthickspace\inf_{x_{n}\in D(T),y_{n+1}\in T(x)}\left\{ \overset{{\scriptscriptstyle \langle n\rangle}}{\pisces}_{T}(x_{n},y)+\langle x-x_{n},y_{n+1}\rangle\right\} .
\end{aligned}
\]
The equalities in (\ref{eq:-60}) and the first equality in (\ref{eq:-61})
follow similar reasoning. Alternatively, in both (\ref{eq:-61}) and
(\ref{eq:-60}), we first prove one equality, and then use it for
$T$ replaced with $T^{-1}$, which then translates into the the other
equality for $T$. Also, (\ref{eq:-61}) implies (\ref{eq:-60}) via
convex conjugation.

If we pass to infimum over $n\ge1$ in (\ref{eq:-61}) we obtain the
corresponding (\ref{eq:-61}) for $n=\infty$; similarly with supremum
for (\ref{eq:-60}). \end{proof}

\strut

In what follows, we apply mathematical induction and the recurrence
relations in (\ref{eq:-61}), (\ref{eq:-60}) to derive closed forms
for the $n-$th Fitzpatrick function $\overset{{\scriptscriptstyle \langle n\rangle}}{\varphi}_{T}$
and $\overset{{\scriptscriptstyle \langle n\rangle}}{\pisces}_{T}$.
These closed forms are then employed to verify or refute (maximal)
$n-$monotonicity with $\overset{{\scriptscriptstyle \langle n\rangle}}{\pisces}_{T}$
being more convenient for establishing monotonicity and $\overset{{\scriptscriptstyle \langle n\rangle}}{\varphi}_{T}$
for demonstrating maximality.

\begin{example} (\textsc{Normal Cone}) Let $(X,Y,\langle\cdot,\cdot\rangle)$
be a dual system and let $C\subset X$. Then 
\begin{equation}
\overset{{\scriptscriptstyle \langle n\rangle}}{\varphi}_{N_{C}}(x,y)=\iota_{C^{\#}}(x)+\sigma_{C}(y),\ \overset{{\scriptscriptstyle \langle n\rangle}}{\psi}_{N_{C}}(x,y)=\iota_{\operatorname*{cl}\operatorname*{conv}C}(x)+\sigma_{C^{\#}}(y),\ x\in X,\ y\in Y,\ n\in\overline{1,\infty};\label{eq:-59}
\end{equation}
\begin{equation}
\begin{aligned}\overset{{\scriptscriptstyle \langle n\rangle}}{\pisces}_{N_{C}}(x,y) & =\iota_{C}(x)+\sigma_{C}(y)+\iota_{R(N_{C})}(y)\\
 & =\left\{ \begin{array}{cc}
\iota_{C}(x)+\sigma_{C}(y), & y\in R(N_{C})\\
+\infty, & \mathit{otherwise}
\end{array}\right.,\ x\in X,\ y\in Y,\ n\in\overline{2,\infty};
\end{aligned}
\label{eq:-64}
\end{equation}
where, for $A\subset X$, $\sigma_{A}(y):=\sup_{u\in A}\langle u,y\rangle$,
and 
\[
C^{\#}:=\{x\in X\mid\forall(u,v)\in\operatorname*{Graph}N_{C},\ \langle x-u,v\rangle\le0\}
\]
is the portable hull of $C$ which is the intersection of all the
supporting half-spaces that contain $C$ and are supported at points
in $C$ (see \cite{Voisei2019}).

As a consequence, $N_{C}$ is maximal cyclically monotone iff $C=C^{\#}$
iff $N_{C}$ is maximal monotone. \end{example}

\begin{proof} If $C=\emptyset$ then $\overset{{\scriptscriptstyle \langle n\rangle}}{\varphi}_{N_{C}}=\sigma_{C}=-\infty$,
$C^{\#}=X$, $\iota_{C^{\#}}=0$, $\overset{{\scriptscriptstyle \langle n\rangle}}{\psi}_{N_{C}}=\iota_{\operatorname*{cl}\operatorname*{conv}C}=\overset{{\scriptscriptstyle \langle n\rangle}}{\pisces}_{N_{C}}=+\infty$,
$\sigma_{C^{\#}}=\iota_{\{0\}}$; whence (\ref{eq:-59}), (\ref{eq:-64})
hold due to the convention $\infty-\infty=\infty$. Assume that $C\neq\emptyset$. 

For $n\in\mathbb{N}$, $n\ge2$, and $(x,y)\in X\times Y$ 
\[
\begin{aligned}\overset{{\scriptscriptstyle \langle n\rangle}}{\varphi}(x,y) & =\sup_{\{(x_{i},y_{i})\}_{i\in\overline{1,n}}\subset\operatorname*{Graph}N_{C}}\{\langle x-x_{1},y_{1}\rangle+\sum_{i=2}^{n}\langle x_{i-1}-x_{i},y_{i}\rangle+\langle x_{n},y\rangle\}\\
 & \le\sup_{(x_{1},y_{1})\in\operatorname*{Graph}N_{C}}\{\langle x-x_{1},y_{1}\rangle\}+\sup_{x_{n}\in C}\langle x_{n},y\rangle=\iota_{C^{\#}}(x)+\sigma_{C}(y);
\end{aligned}
\]
the inequality is verified similarly for $n=1$ and, as a consequence,
holds for $n=\infty$. 

For every $n\in\overline{1,\infty}$, $(x,y)\in\operatorname*{dom}\overset{{\scriptscriptstyle \langle n\rangle}}{\varphi}_{N_{C}}\subset\operatorname*{dom}\overset{{\scriptscriptstyle \langle1\rangle}}{\varphi}_{N_{C}}$
we have 
\[
\forall a\in C,\ n^{*}\in N_{C}(a),\ \langle x-a,n^{*}\rangle+\langle a,y\rangle\le\overset{{\scriptscriptstyle \langle1\rangle}}{\varphi}_{N_{C}}(x,y)<+\infty.
\]
Since, for every $a\in C,$ $N_{C}(a)$ is a cone, we find that $x\in C^{\#}$.
Because, for every $a\in C,$ $0\in N_{C}(a)$ we infer that $\overset{{\scriptscriptstyle \langle n\rangle}}{\varphi}_{N_{C}}(x,y)\ge\overset{{\scriptscriptstyle \langle1\rangle}}{\varphi}_{N_{C}}(x,y)\ge\sigma_{C}(y)$.
Hence $\overset{{\scriptscriptstyle \langle n\rangle}}{\varphi}_{N_{C}}(x,y)\ge\iota_{C^{\#}}(x)+\sigma_{C}(y)$.
Therefore the first part in (\ref{eq:-59}) holds. The second part
in (\ref{eq:-59}) follows by convex conjugation.

To prove (\ref{eq:-64}) we use induction over $n\ge2$. According
to (\ref{eq:-5}), 
\[
\overset{{\scriptscriptstyle \langle2\rangle}}{\pisces}_{N_{C}}(x,y)=\inf_{x_{1}\in(N_{C})^{-1}(y),\ y_{2}\in N_{C}(x)}\{\langle x_{1},y\rangle+\langle x-x_{1},y_{2}\rangle\}.
\]
To complete the basis step, it suffices to prove $\overset{{\scriptscriptstyle \langle2\rangle}}{\pisces}_{N_{C}}(x,y)=\sigma_{C}(y)$
when $x\in C$, $y\in R(N_{C})$. To this end note that, when $x\in C$,
for every $y_{2}\in N_{C}(x)$, $\langle x-x_{1},y_{2}\rangle\ge0$,
because $x_{1}\in C$. Therefore, since $0\in N_{C}(x)$, $\inf_{y_{2}\in N_{C}(x)}\langle x-x_{1},y_{2}\rangle=0$
and so $\overset{{\scriptscriptstyle \langle2\rangle}}{\pisces}_{N_{C}}(x,y)=\inf_{x_{1}\in(N_{C})^{-1}(y)}\langle x_{1},y\rangle=\sigma_{C}(y)$.

Assume that $\overset{{\scriptscriptstyle \langle n\rangle}}{\pisces}_{N_{C}}(x,y)=\iota_{C}(x)+\sigma_{C}(y)+\iota_{R(N_{C})}(y)$,
for some integer $n\ge2$. From \ref{eq:-61} 
\[
\overset{{\scriptscriptstyle \langle n+1\rangle}}{\pisces}_{N_{C}}(x,y)=\sigma_{C}(y)+\iota_{R(N_{C})}(y)+\inf_{u\in C,v\in N_{C}(x)}\langle x-u,v\rangle=\iota_{C}(x)+\sigma_{C}(y)+\iota_{R(N_{C})}(y),
\]
since, for $x\in C$, $\inf_{u\in C,v\in N_{C}(x)}\left\{ \langle x-u,v\rangle\right\} =0$.
The inductive step is complete. 

Hence (\ref{eq:-64}) holds for every $n\in\mathbb{N}$, $n\ge2$
and also for $n=\infty$.

Assume that $N_{C}$ is maximal cyclically monotone. Since $\operatorname*{Graph}(N_{C})\subset\operatorname*{Graph}(N_{\operatorname*{cl}\operatorname*{conv}C})$,
we know that $C$ is closed convex. According to Theorem \ref{F},
\[
\iota_{C^{\#}}=\iota_{C}^{\cup}=\iota_{C}^{\circ}=\iota_{C},
\]
that is, $C^{\#}=C$. Hence, for every $x\in X$, $y\in Y$, we have
$\overset{{\scriptscriptstyle \langle1\rangle}}{\varphi}_{N_{C}}(x,y)=\iota_{C}(x)+\sigma_{C}(y)\ge\langle x,y\rangle$.
We conclude that $N_{C}$ is maximal monotone with the aid of Theorem
\ref{fzv} (ii). \end{proof}

\begin{lemma} Let $(H,\cdot)$ be a Hilbert space with induced norm
$\|\cdot\|$. For every $v\in H$, 
\begin{equation}
\sup_{u\in H}\{\alpha\|u\|^{2}\pm u\cdot v\}=\left\{ \begin{array}{cc}
-\tfrac{1}{4\alpha}\|v\|^{2}, & \alpha<0\\
\iota_{\{0\}}(v), & \alpha=0\\
\infty, & \alpha>0
\end{array}\right.,\ \inf_{u\in H}\{\alpha\|u\|^{2}\pm u\cdot v\}=\left\{ \begin{array}{cc}
-\tfrac{1}{4\alpha}\|v\|^{2}, & \alpha>0\\
-\iota_{\{0\}}(v), & \alpha=0\\
-\infty, & \alpha<0
\end{array}\right..\label{eq:-69}
\end{equation}
\end{lemma}

\begin{proof} Note that, for $a\neq0,$ $u,v\in H$, we have $a\|u\|^{2}\pm u\cdot v=a\|u\pm\frac{1}{2a}v\|^{2}-\tfrac{1}{4a}\|v\|^{2}$.
\end{proof}

\begin{example} \label{Id} (\textsc{Identity}) Let $(H,\cdot)$
be a Hilbert space with induced norm $\|\cdot\|$ and let $I:H\to H$,
$I(x)=x$. For every $x,y\in H$ 
\begin{equation}
\overset{{\scriptscriptstyle \langle n\rangle}}{\varphi}_{I}(x,y)=\tfrac{n}{2(n+1)}\|x-y\|^{2}+x\cdot y,\ n\ge1;\label{eq:-62}
\end{equation}
\begin{equation}
\overset{{\scriptscriptstyle \langle1\rangle}}{\psi}_{I}(x,y)=\overset{{\scriptscriptstyle \langle1\rangle}}{\pisces}_{I}=c_{I},\ \overset{{\scriptscriptstyle \langle n\rangle}}{\psi}_{I}(x,y)=\overset{{\scriptscriptstyle \langle n\rangle}}{\pisces}_{I}(x,y)=\tfrac{n}{2(n-1)}\|x-y\|^{2}+x\cdot y,\ n\ge2;\label{eq:-63}
\end{equation}
\begin{equation}
\overset{{\scriptscriptstyle \langle\infty\rangle}}{\varphi}_{I}(x,y)=\overset{{\scriptscriptstyle \langle\infty\rangle}}{\pisces}_{I}(x,y)=\tfrac{1}{2}\|x\|^{2}+\tfrac{1}{2}\|y\|^{2};\label{eq:-68}
\end{equation}
and $I$ is cyclically and maximal monotone. \end{example}

\begin{proof} For (\ref{eq:-62}) we use induction on $n\ge1$. From
(\ref{eq:-69}), for every $x,y\in H$ 
\[
\overset{{\scriptscriptstyle \langle1\rangle}}{\varphi}_{I}(x,y)=\sup_{u\in H}\{(x-u)\cdot u+u\cdot y\}=\sup_{u\in H}\{-\|u\|^{2}+u\cdot(x+y)\}=\tfrac{1}{4}\|x+y\|^{2}=\tfrac{1}{4}\|x-y\|^{2}+x\cdot y.
\]
Assume that, for every $x,y\in H$, $\overset{{\scriptscriptstyle \langle n\rangle}}{\varphi}_{I}(x,y)=\frac{n}{2(n+1)}\|x-y\|^{2}+x\cdot y$,
for some finite $n\ge1$. From \ref{eq:-60}, (\ref{eq:-69}) we get
\[
\begin{aligned}\overset{{\scriptscriptstyle \langle n+1\rangle}}{\varphi}_{I}(x,y) & =\sup_{w\in H}\left\{ \tfrac{n}{2(n+1)}\|w-y\|^{2}+w\cdot y+(x-w)\cdot w\right\} \\
=\sup_{u\in H} & \left\{ \tfrac{n}{2(n+1)}\|u\|^{2}+(y+u)\cdot y+(x-y-u)\cdot(y+u)\right\} \\
=\sup_{u\in H} & \left\{ -\tfrac{n+2}{2(n+1)}\|u\|^{2}+u\cdot(x-y)\right\} +x\cdot y=\tfrac{n+1}{2(n+2)}\|x_{}-y\|^{2}+x\cdot y.
\end{aligned}
\]
This completes the inductive step and proves (\ref{eq:-62}).

We proceed similarly to prove, by induction on $n\ge2$, that 
\begin{equation}
\overset{{\scriptscriptstyle \langle n\rangle}}{\pisces}_{I}(x,y)=\tfrac{n}{2(n-1)}\|x-y\|^{2}+x\cdot y.\label{eq:-70}
\end{equation}
By a direct computation $\overset{{\scriptscriptstyle \langle2\rangle}}{\pisces}_{I}(x,y)=\|x\|^{2}+\|y\|^{2}-x\cdot y=\|x-y\|^{2}+x\cdot y$.

Assume that (\ref{eq:-70}) is true for some finite integer $n\ge2$.
Use (\ref{eq:-61}), (\ref{eq:-69}) to get 
\[
\begin{aligned}\overset{{\scriptscriptstyle \langle n+1\rangle}}{\pisces}_{I}(x,y) & =\inf_{w\in H}\left\{ \tfrac{n}{2(n-1)}\|w-y\|^{2}+w\cdot(y-x)+\|x\|^{2}\right\} \\
 & =\inf_{u\in H}\left\{ \tfrac{n}{2(n-1)}\|u\|^{2}+u\cdot(y-x)\right\} +\|x-y\|^{2}+x\cdot y\\
 & =-\tfrac{n-1}{2n}\|x-y\|^{2}+\|x-y\|^{2}+x\cdot y=\tfrac{n+1}{2n}\|x-y\|^{2}+x\cdot y.
\end{aligned}
\]

Note that $\overset{{\scriptscriptstyle \langle n\rangle}}{\pisces}_{I}\in\Gamma(H)$,
e.g. because its Gâteaux derivative $\nabla\overset{{\scriptscriptstyle \langle n\rangle}}{\pisces}_{I}(x,y)=\frac{1}{(n-1)}(nx-y,ny-x)$,
$x,y\in H$, is linear and monotone with 
\[
\langle\nabla\overset{{\scriptscriptstyle \langle n\rangle}}{\pisces}_{I}(x,y),(x,y)\rangle_{H\times H}=\|x\|^{2}+\|y\|^{2}+\tfrac{1}{n-1}\|x-y\|^{2}\ge0,\ x,y\in H,
\]

Hence $\overset{{\scriptscriptstyle \langle n\rangle}}{\psi}_{I}=\overset{{\scriptscriptstyle \langle n\rangle}}{\pisces}_{I}$
and the proof of (\ref{eq:-63}) is complete. \end{proof}

\begin{example} (\textsc{Rotation}) Let $R_{\theta}$ be the counterclockwise
vector rotation in $\mathbb{R}^{2}$ by an angle of $\theta\in[-\pi,\pi]$,
i.e., $R_{\theta}:\mathbb{R}^{2}\to\mathbb{R}^{2}$, 
\[
R_{\theta}(x)=\begin{pmatrix}\cos(\theta) & -\sin(\theta)\\
\sin(\theta) & \cos(\theta)
\end{pmatrix}\begin{pmatrix}x_{1}\\
x_{2}
\end{pmatrix}=\begin{pmatrix}x_{1}\cos\theta-x_{2}\sin\theta\\
x_{1}\sin\theta+x_{2}\cos\theta
\end{pmatrix},\ x=\begin{pmatrix}x_{1}\\
x_{2}
\end{pmatrix}.
\]

For $x,y\in\mathbb{R}^{2}$ 
\begin{equation}
\overset{{\scriptscriptstyle \langle n\rangle}}{\varphi}_{R_{\theta}}(x,y)=\left\{ \begin{array}{cc}
\frac{n}{2(n+1)}\|y-R_{\theta}(x)\|^{2}+x\cdot y, & \theta=0\\
\frac{\sin(n\theta)}{2\sin((n+1)\theta)}\|y-R_{\theta}(x)\|^{2}+x\cdot y, & \theta\in(-\frac{\pi}{n+1},\frac{\pi}{n+1})\smallsetminus\{0\}\\
c_{R_{\theta}}(x,y), & \theta=\pm\frac{\pi}{n+1}\\
+\infty, & \theta\in[-\pi,-\frac{\pi}{n+1})\cup(\frac{\pi}{n+1},\pi]
\end{array}\right.,\ n\in\mathbb{N},n\ge1;\label{eq:-65}
\end{equation}
\begin{equation}
\overset{{\scriptscriptstyle \langle\infty\rangle}}{\varphi}_{R_{\theta}}(x,y)=\left\{ \begin{array}{cc}
\frac{1}{2}\|x\|^{2}+\frac{1}{2}\|y\|^{2}, & \theta=0\\
+\infty, & \theta\in[-\pi,\pi]\smallsetminus\{0\}
\end{array}\right.;\label{eq:-73}
\end{equation}
\begin{equation}
\overset{{\scriptscriptstyle \langle2\rangle}}{\pisces}_{R_{\theta}}(x,y)=\cos\theta\|y-R_{\theta}(x)\|^{2}+x\cdot y,\ \theta\in[-\pi,\pi];\label{eq:-71}
\end{equation}
\begin{equation}
\overset{{\scriptscriptstyle \langle n\rangle}}{\pisces}_{R_{\theta}}(x,y)=\left\{ \begin{array}{cc}
\frac{n}{2(n-1)}\|y-R_{\theta}(x)\|^{2}+x\cdot y, & \theta=0\\
\frac{\sin(n\theta)}{2\sin((n-1)\theta)}\|y-R_{\theta}(x)\|^{2}+x\cdot y, & \theta\in(-\frac{\pi}{n-1},\frac{\pi}{n-1})\smallsetminus\{0\}\\
c-\iota_{\{y-R_{\theta}(x)\}}, & \theta=\pm\frac{\pi}{n-1}\\
-\infty, & \theta\in[-\pi,-\frac{\pi}{n-1})\cup(\frac{\pi}{n-1},\pi]
\end{array}\right.,\ n\in\mathbb{N},n\ge3;\label{eq:-72}
\end{equation}
\begin{equation}
\overset{{\scriptscriptstyle \langle\infty\rangle}}{\pisces}_{R_{\theta}}(x,y)=\left\{ \begin{array}{cc}
\frac{1}{2}\|x\|^{2}+\frac{1}{2}\|y\|^{2}, & \theta=0\\
-\infty, & \theta\in[-\pi,\pi]\smallsetminus\{0\}
\end{array}\right..\label{eq:-66}
\end{equation}
For $n\in\mathbb{N}$, $n\ge2$, $R_{\theta}$ is (maximal) $n-$monotone
iff $\theta\in[-\frac{\pi}{n},\frac{\pi}{n}]$; $R_{\theta}$ is $n-$monotone
and not $(n+1)-$monotone iff $\theta\in[-\frac{\pi}{n},-\frac{\pi}{n+1})\cup(\frac{\pi}{n+1},\frac{\pi}{n}]$;
and $R_{\theta}$ is (maximal) cyclically monotone iff $\theta=0$.
\end{example}

\begin{proof} We repeatedly use the fact that both the norm and inner
product in $\mathbb{R}^{2}$ are invariant under any rotation, namely,
$\|R_{\theta}(x)\|=\|x\|$ and $R_{\theta}(x)\cdot R_{\theta}(y)=x\cdot y$,
$x,y\in\mathbb{R}^{2}$.

By a direct computation $x\cdot R_{\theta}(x)=\cos\theta\|x\|^{2}$
and so, together with $R_{\theta}^{*}=R_{\theta}^{-1}=R_{-\theta}$,
\[
\begin{aligned}\overset{{\scriptscriptstyle \langle2\rangle}}{\pisces}_{R_{\theta}}(x,y) & =R_{-\theta}(y)\cdot y+(x-R_{-\theta}(y))\cdot R_{\theta}(x)\\
 & =(y-R_{\theta}(x))\cdot R_{\theta}(y-R_{\theta}(x))+x\cdot y=\cos\theta\|y-R_{\theta}(x)\|^{2}+x\cdot y
\end{aligned}
\]

Therefore, according to Theorem \ref{cm0}, $R_{\theta}$ is ($2-$)monotone
iff $\theta\in[-\frac{\pi}{2},\frac{\pi}{2}]$; in this case, it is
maximal ($2-$)monotone because it is monotone and defined on $\mathbb{R}^{2}$. 

Alternatively, $R_{\theta}$ is maximal monotone for $\theta\in[-\frac{\pi}{2},\frac{\pi}{2}]$,
since $c_{R_{\theta}}\in\Gamma(\mathbb{R}^{2})$, which implies that
$\overset{{\scriptscriptstyle \langle1\rangle}}{\psi}_{R_{\theta}}=\overset{{\scriptscriptstyle \langle1\rangle}}{\pisces}_{R_{\theta}}=c_{R_{\theta}}\ge c$.
Hence $R_{\theta}$ is $2-$representable, due to $\operatorname*{Graph}R_{\theta}=[\overset{{\scriptscriptstyle \langle1\rangle}}{\psi}_{R_{\theta}}=c]$
and $\overset{{\scriptscriptstyle \langle1\rangle}}{\psi}_{R_{\theta}}\ge c$
(see Theorem \ref{r2inf}). In addition, from (\ref{eq:-65}), $\overset{{\scriptscriptstyle \langle1\rangle}}{\varphi}_{R_{\theta}}\ge c$.
Theorem \ref{fzv} completes the argument.

We restrict ourselves to $\theta\in(0,\pi]$, since $R_{0}=I$ is
covered in Example \ref{Id} and we can extend our conclusions to
$\theta\in[-\pi,0)$ via $R_{-\theta}=R_{\theta}^{-1}=R_{\theta}^{*}$
and $\overset{{\scriptscriptstyle \langle n\rangle}}{\varphi}_{R_{\theta}}(x,y)=\overset{{\scriptscriptstyle \langle n\rangle}}{\varphi}_{R_{-\theta}}(y,x)$,
$\overset{{\scriptscriptstyle \langle n\rangle}}{\pisces}_{R_{\theta}}(x,y)=\overset{{\scriptscriptstyle \langle n\rangle}}{\pisces}_{R_{-\theta}}(y,x)$.

According to (\ref{eq:-61}), (\ref{eq:-69}) 
\[
\begin{aligned}\overset{{\scriptscriptstyle \langle3\rangle}}{\pisces}_{R_{\theta}}(x,y) & =\inf_{v\in\mathbb{R}^{2}}\left\{ \overset{{\scriptscriptstyle \langle2\rangle}}{\pisces}_{R_{\theta}}(x,v)+R_{\theta}^{-1}(y)\cdot(y-v)\right\} \\
 & =\inf_{v\in\mathbb{R}^{2}}\left\{ \cos\theta\|v-R_{\theta}(x)\|^{2}+v\cdot(x-R_{\theta}^{-1}(y))+y\cdot R_{\theta}(y)\right\} \\
 & =\inf_{u\in\mathbb{R}^{2}}\left\{ \cos\theta\|u\|^{2}+u\cdot(x-R_{\theta}^{-1}(y))\right\} +\overset{{\scriptscriptstyle \langle2\rangle}}{\pisces}_{R_{\theta}}(x,y)\\
 & \inf_{u\in\mathbb{R}^{2}}\left\{ \cos\theta\|u\|^{2}+u\cdot(x-R_{\theta}^{-1}(y))\right\} +\cos\theta\|y-R_{\theta}(x)\|^{2}+x\cdot y\\
 & =\left\{ \begin{array}{cc}
(\cos\theta-\frac{1}{4\cos\theta})\|y-R_{\theta}(x)\|^{2}+x\cdot y, & \cos\theta>0\\
c-\iota_{\{y-R_{\theta}(x)\}}, & \cos\theta=0\\
-\infty & \cos\theta<0
\end{array}\right.\\
 & =\left\{ \begin{array}{cc}
\frac{\sin(3\theta)}{2\sin(2\theta)}\|y-R_{\theta}(x)\|^{2}+x\cdot y, & \theta\in(0,\frac{\pi}{2})\\
c-\iota_{\{y-R_{\theta}(x)\}}, & \theta=\frac{\pi}{2}\\
-\infty & \theta\in(\frac{\pi}{2},\pi]
\end{array},\right.
\end{aligned}
\]
after we used the identity $\cos\theta-\frac{1}{4\cos\theta}=\frac{\sin(3\theta)}{2\sin(2\theta)}$.
Hence the base case for an induction argument on $n\ge3$ used for
(\ref{eq:-72}) is complete.

Assume that (\ref{eq:-72}) holds for some finite integer $n\ge3$.
From (\ref{eq:-61}) 
\[
\overset{{\scriptscriptstyle \langle n+1\rangle}}{\pisces}_{R_{\theta}}(x,y)=\inf_{v\in\mathbb{R}^{2}}\left\{ \overset{{\scriptscriptstyle \langle n\rangle}}{\pisces}_{T}(x,v)+R_{\theta}^{-1}(y)\cdot(y-v)\right\} ;
\]
from which $\overset{{\scriptscriptstyle \langle n+1\rangle}}{\pisces}_{R_{\theta}}(x,y)=-\infty$,
whenever $\theta\in[\frac{\pi}{n-1},\pi]$.

Therefore, for $\theta\in(0,\frac{\pi}{n-1})$ 
\[
\begin{aligned}\overset{{\scriptscriptstyle \langle n+1\rangle}}{\pisces}_{R_{\theta}}(x,y) & =\inf_{v\in\mathbb{R}^{2}}\left\{ \tfrac{\sin(n\theta)}{2\sin((n-1)\theta)}\|v-R_{\theta}(x)\|^{2}+v\cdot(x-R_{\theta}^{-1}(y))+y\cdot R_{\theta}(y)\right\} \\
 & =\inf_{u\in\mathbb{R}^{2}}\left\{ \tfrac{\sin(n\theta)}{2\sin((n-1)\theta)}\|u\|^{2}+u\cdot(x-R_{\theta}^{-1}(y))\right\} +\overset{{\scriptscriptstyle \langle2\rangle}}{\pisces}_{R_{\theta}}(x,y)\\
 & =\left\{ \begin{array}{cc}
(\cos\theta-\frac{\sin((n-1)\theta)}{2\sin(n\theta)})\|y-R_{\theta}(x)\|^{2}+x\cdot y, & \frac{\sin(n\theta)}{\sin((n-1)\theta)}>0\\
c-\iota_{\{y-R_{\theta}(x)\}}, & \sin(n\theta)=0\\
-\infty & \frac{\sin(n\theta)}{\sin((n-1)\theta)}<0
\end{array}\right.\\
 & =\left\{ \begin{array}{cc}
\frac{\sin((n+1)\theta)}{2\sin(n\theta)}\|y-R_{\theta}(x)\|^{2}+x\cdot y, & \theta\in(0,\frac{\pi}{n})\\
c-\iota_{\{y-R_{\theta}(x)\}}, & \theta=\frac{\pi}{n}\\
-\infty & \theta\in(\frac{\pi}{n},\frac{\pi}{n-1})
\end{array}\right.
\end{aligned}
\]
where the identity $\cos\theta-\frac{\sin((n-1)\theta)}{2\sin(n\theta)}=\frac{\sin((n+1)\theta)}{2\sin(n\theta)}$
is used and, since $\theta\in(0,\frac{\pi}{n-1})$, i.e., $\sin((n-1)\theta)>0$,
we know that the sign of $\frac{\sin(n\theta)}{\sin((n-1)\theta)}$
is given by the sign of $\sin(n\theta)$.

Hence (\ref{eq:-72}) holds for every finite integer $n\ge3$.

We establish (\ref{eq:-65}) by applying a similar argument. For $n=1$
\[
\begin{aligned}\overset{{\scriptscriptstyle \langle1\rangle}}{\varphi}_{R_{\theta}}(x,y) & =\sup_{u\in\mathbb{R}^{2}}\{(x-u)\cdot R_{\theta}(u)+u\cdot y\}=\sup_{u\in\mathbb{R}^{2}}\{-\cos\theta\|u\|^{2}+u\cdot(y+R_{\theta}^{*}(x))\}\\
 & =\left\{ \begin{array}{cc}
\tfrac{1}{4\cos\theta}\|y+R_{\theta}^{*}(x)\|^{2}, & \cos\theta>0\\
\iota_{\{0\}}(y+R_{\theta}^{*}(x)), & \cos\theta=0\\
+\infty, & \cos\theta<0
\end{array}\right.=\left\{ \begin{array}{cc}
\tfrac{1}{4\cos\theta}\|y-R_{\theta}(x)\|^{2}+x\cdot y, & \theta\in(0,\frac{\pi}{2})\\
c_{R_{\frac{\pi}{2}}}(x,y), & \theta=\frac{\pi}{2}\\
+\infty, & \theta\in(\frac{\pi}{2},\pi]
\end{array}\right.
\end{aligned}
\]
since $R_{\theta}(x)+R_{\theta}^{*}(x)=(2\cos\theta)x$, $R_{\frac{\pi}{2}}^{*}=R_{-\frac{\pi}{2}}=-R_{\frac{\pi}{2}}$,
and $\tfrac{1}{4\cos\theta}=\tfrac{\sin\theta}{2\sin(2\theta)}$.

Assume that (\ref{eq:-65}) holds for a finite integer $n\ge1$. From
(\ref{eq:-61}) 
\[
\overset{{\scriptscriptstyle \langle n+1\rangle}}{\varphi}_{R_{\theta}}(x,y)=\sup_{u\in\mathbb{R}^{2}}\left\{ \overset{{\scriptscriptstyle \langle n\rangle}}{\varphi}_{R_{\theta}}(u,y)+(x-u)\cdot R_{\theta}(u)\right\} ;
\]
from which $\overset{{\scriptscriptstyle \langle n+1\rangle}}{\varphi}_{R_{\theta}}(x,y)=+\infty$,
whenever $\theta\in[\frac{\pi}{n+1},\pi]$.

For $\theta\in(0,\frac{\pi}{n+1})$, $\sin((n+1)\theta)>0$. According
to (\ref{eq:-61}), (\ref{eq:-69}) and the identity 
\[
\tfrac{\sin(n\theta)}{2\sin((n+1)\theta)}-\cos\theta=-\tfrac{\sin((n+2)\theta)}{2\sin((n+1)\theta)}
\]
we have 
\[
\begin{aligned}\overset{{\scriptscriptstyle \langle n+1\rangle}}{\varphi}_{R_{\theta}}(x,y) & =\sup_{u\in\mathbb{R}^{2}}\left\{ \tfrac{\sin(n\theta)}{2\sin((n+1)\theta)}\|R_{\theta}(u)-R_{\theta}(x)\|^{2}+x\cdot R_{\theta}(u)+u\cdot(y-R_{\theta}(u))\right\} \\
 & =\sup_{u\in\mathbb{R}^{2}}\left\{ \tfrac{\sin(n\theta)}{2\sin((n+1)\theta)}\|u-x\|^{2}+x\cdot R_{\theta}(u)+u\cdot(y-R_{\theta}(u))\right\} \\
 & =\sup_{v\in\mathbb{R}^{2}}\left\{ \tfrac{\sin(n\theta)}{2\sin((n+1)\theta)}\|v\|^{2}-v\cdot R_{\theta}(v)+v\cdot(y-R_{\theta}(x))+x\cdot y\right\} \\
 & =\sup_{v\in\mathbb{R}^{2}}\left\{ (\tfrac{\sin(n\theta)}{2\sin((n+1)\theta)}-\cos\theta)\|v\|^{2}+v\cdot(y-R_{\theta}(x))\right\} +x\cdot y\\
 & =\sup_{v\in\mathbb{R}^{2}}\left\{ -\tfrac{\sin((n+2)\theta)}{2\sin((n+1)\theta)}\|v\|^{2}+v\cdot(y-R_{\theta}(x))\right\} +x\cdot y\\
 & =\left\{ \begin{array}{cc}
\tfrac{\sin((n+1)\theta)}{2\sin((n+2)\theta)}\|y-R_{\theta}(x)\|^{2}+x\cdot y, & \sin((n+2)\theta)>0\\
\iota_{\{0\}}(y-R_{\theta}(x))+x\cdot y, & \sin((n+2)\theta)=0\\
+\infty & \sin((n+2)\theta)<0
\end{array}\right.\\
 & =\left\{ \begin{array}{cc}
\tfrac{\sin((n+1)\theta)}{2\sin((n+2)\theta)}\|y-R_{\theta}(x)\|^{2}+x\cdot y, & \theta\in(0,\frac{\pi}{n+2})\\
c_{R_{\theta}}(x,y), & \theta=\frac{\pi}{n+2}\\
+\infty & \theta\in(\frac{\pi}{n+2},\frac{\pi}{n+1})
\end{array}\right.
\end{aligned}
\]
and the inductive step is complete. Therefore (\ref{eq:-65}) is true.

Note that, from (\ref{eq:-72}), for $n\in\mathbb{N}$, $n\ge3$,
$R_{\theta}$ is $n-$monotone iff $\overset{{\scriptscriptstyle \langle n\rangle}}{\pisces}_{R_{\theta}}\ge c$
which comes to $\theta\in(-\frac{\pi}{n-1},\frac{\pi}{n-1})$ and
either $\theta=0$ or $\sin(n\theta)>0$, that is, $\theta\in[-\frac{\pi}{n},\frac{\pi}{n}]$.
\end{proof}

\strut

A subset $A$ of $Z=X\times Y$ is a \emph{double-cone} if $\mathbb{R}A=A$,
i.e., for every $\lambda\in\mathbb{R}$, $z\in A$ we have $\lambda z\in A$.
A subset $A$ of $Z$ is called \emph{non-negative} if, for every
$z\in A$, $c(z)\geq0$, and \emph{skew} if, for every $z\in A$,
$c(z)=0$.

For $A\subset Z$ we set 
\begin{align*}
\neg A & :=\{z=(x,y)\in Z\mid\neg z:=(x,-y)\in A\},\\
A^{\perp} & :=\{z^{\prime}\in Z\mid\forall z\in A,\ z\cdot z^{\prime}=0\},\ \mathrm{and}
\end{align*}
Note that $\neg A\neq-A:=(-1)A$, $A^{\perp}=(\operatorname*{cl}_{\sigma(Z,Z)}\operatorname*{lin}\!A)^{\perp}$,
$A^{\perp\perp}=\operatorname*{cl}_{\sigma(Z,Z)}\operatorname*{lin}\!A$,
$A^{\perp}$ is a $\sigma(Z,Z)-$closed linear subspace of $Z$, and
$\neg(A^{\perp})=(\neg A)^{\perp}$.

\begin{example} \label{SDC} (\textsc{Skew Double-Cone}) Let $(X,Y,\langle\cdot,\cdot\rangle)$
be a dual system and let $D\subset Z:=X\times Y$ be a skew double-cone.
Then

\emph{(i)} 
\begin{equation}
\overset{{\scriptscriptstyle \langle1\rangle}}{\pisces}_{D}=\iota_{D},\ \overset{{\scriptscriptstyle \langle1\rangle}}{\varphi}_{D}=\iota_{D^{\perp}},\ \overset{{\scriptscriptstyle \langle1\rangle}}{\psi}_{D}=\iota_{D^{\perp\perp}},\label{eq:-77}
\end{equation}
$D$ is monotone iff $D\subset D^{\perp}$; $D$ is representable
iff $D=D^{\perp\perp}$ iff $D$ is a $\sigma(Z,Z)-$closed linear
subspace of $Z$; $D$ is maximal monotone iff $D$ is a $\sigma(Z,Z)-$closed
linear subspace of $Z$ and $\neg D^{\perp}$ is monotone.

\medskip{}

\emph{(ii)} 
\begin{equation}
\overset{{\scriptscriptstyle \langle2\rangle}}{\varphi}_{D}(x,y)=\sup_{u\in\Pr_{X}(D)}\iota_{D^{\perp}}(u,y)+\iota_{(\Pr_{Y}(D))^{\perp}}(x),\ (x,y)\in Z;\label{eq:-78}
\end{equation}
where $\Pr_{X}(x,y):=x$ and $\Pr_{Y}(x,y):=y$, $(x,y)\in X\times Y$
are the projections of $X\times Y$ onto $X$ and $Y$, respectively.

\medskip{}

\emph{(iii)} $D$ is $3-$monotone iff $\Pr_{X}(D)\subset(\Pr_{Y}(D))^{\perp}$
iff, for every $(x,y)\in\Pr_{X}(D)\times\Pr_{Y}(D)$, $\langle x,y\rangle=0$
iff $D\subset S\times S^{\perp}$, for some $S$ a $\sigma(X,Y)-$closed
linear subspace of $X$ iff $D$ is cyclically monotone.

\medskip{}

\emph{(iv)} $D$ is maximal $3-$monotone iff $D=S\times S^{\perp}$,
for some $S$ a $\sigma(X,Y)-$closed linear subspace of $X$ iff
$D$ is maximal cyclically monotone. In this case $\overset{{\scriptscriptstyle \langle2\rangle}}{\varphi}_{D}=\overset{{\scriptscriptstyle \langle\infty\rangle}}{\varphi}_{D}=\iota_{S\times S^{\perp}}$.
\end{example}

\begin{proof} (i) Since $D$ is a skew double-cone, we have $\overset{{\scriptscriptstyle \langle1\rangle}}{\pisces}_{D}=c_{D}=\iota_{D}$,
and 
\[
\overset{{\scriptscriptstyle \langle1\rangle}}{\varphi}_{D}(z)=\sup\{z\cdot z^{\prime}-c(z^{\prime})\mid z^{\prime}\in D\}=\sup\{z\cdot z^{\prime}\mid z^{\prime}\in D\}=\iota_{D^{\perp}}(z),\ z\in Z.
\]
Because $D^{\perp}$ is a linear subspace of $Z$ we get $\overset{{\scriptscriptstyle \langle1\rangle}}{\psi}_{D}=(\overset{{\scriptscriptstyle \langle1\rangle}}{\varphi}_{D})^{\square}=\iota_{D^{\perp}}^{\square}=\sigma_{D^{\perp}}=\iota_{D^{\perp\perp}}$.

According to Theorem \ref{cm1}, $D$ is monotone iff $\iota_{D^{\perp}}=\overset{{\scriptscriptstyle \langle1\rangle}}{\varphi}_{D}\le\overset{{\scriptscriptstyle \langle1\rangle}}{\pisces}_{D}=\iota_{D}$
iff $D\subset D^{\perp}$.

According to Theorem \ref{r2inf}, $D$ is representable iff $\overset{{\scriptscriptstyle \langle1\rangle}}{\psi}_{D}=\iota_{D^{\perp\perp}}\ge c$
and $D=[\iota_{D^{\perp\perp}}=c]$. Clearly, since $D$ is skew,
if $D=D^{\perp\perp}$ then $D$ is representable.

Conversely, assume that $D$ is representable, that is, for every
$z\in L:=D^{\perp\perp}=\operatorname*{cl}_{\sigma(Z,Z)}\operatorname*{lin}D$,
we have $c(z)\le0$, and $D=[\iota_{L}=c]$. It suffices to prove
that $L$ is skew, since, in that case, $[\iota_{L}=c]=L$. To this
end note that, because $D$ is skew monotone (recall (\ref{eq:-20})),
for every $z,w\in D$, we have 
\[
c(z-w)=c(z)+c(w)+z\cdot w=z\cdot w\ge0.
\]
Since $-z\in D$, we know that $z\cdot w=0$. Therefore, because $\neg D$
is skew, for every $u,v\in\neg D$ 
\[
c(u-v)=c(u)+c(v)-u\cdot v=\neg u\cdot\neg v=0,
\]
in particular, $\neg D$ is a skew monotone double-cone. According
to (\ref{eq:-77}) and Theorem \ref{m2}, $\overset{{\scriptscriptstyle \langle1\rangle}}{\psi}_{\neg D}=\iota_{(\neg D)^{\perp\perp}}=\iota_{\neg L}\ge c$,
which comes to, for every $z\in\neg L$, $c(z)\le0$. Hence, for every
$z\in L$, $c(z)=-c(\neg z)\ge0$, and so, $L$ is skew.

According to Theorem \ref{fzv}, $D$ is maximal monotone iff $\overset{{\scriptscriptstyle \langle1\rangle}}{\varphi}_{D}=\iota_{D^{\perp}}\ge c$
and $D$ is representable. Because $D^{\perp}$ is linear, $\iota_{D^{\perp}}\ge c$
is equivalent to $\neg D^{\perp}$ is monotone.

(ii) From (\ref{eq:-60}) and the fact that $D$ is a skew double-cone,
we find that 
\[
\begin{aligned}\overset{{\scriptscriptstyle \langle2\rangle}}{\varphi}_{D}(x,y) & =\sup_{(u,v)\in D}\{\iota_{D^{\perp}}(u,y)+\langle x-u,v\rangle\}=\sup_{(u,v)\in D}\{\iota_{D^{\perp}}(u,y)+\langle x,v\rangle\}\\
 & =\sup_{u\in\Pr_{X}(D)}\iota_{D^{\perp}}(u,y)+\sigma_{\Pr_{Y}(D)}(x)=\sup_{u\in\Pr_{X}(D)}\iota_{D^{\perp}}(u,y)+\iota_{(\Pr_{Y}(D))^{\perp}}(x),\ (x,y)\in Z,
\end{aligned}
\]
because $\Pr_{X}(D)$ is a double-cone.

(iii) According to Theorem \ref{cm1}, $D$ is $3-$monotone iff $\overset{{\scriptscriptstyle \langle2\rangle}}{\varphi}_{D}\le\overset{{\scriptscriptstyle \langle1\rangle}}{\pisces}_{D}=\iota_{D}$
iff, for every $z\in D$, $\overset{{\scriptscriptstyle \langle2\rangle}}{\varphi}_{D}\le0$.

If $\Pr_{X}(D)\subset(\Pr_{Y}(D))^{\perp}$, i.e., for every $(x,y)\in\Pr_{X}(D)\times\Pr_{Y}(D)$,
$\langle x,y\rangle=0$, then, for every $(x,y)\in D$, $u\in\Pr_{X}(D)$,
we have $(u,y)\in D^{\perp}$ and $x\in\Pr_{X}(D)\subset(\Pr_{Y}(D))^{\perp}$.
Hence $\overset{{\scriptscriptstyle \langle2\rangle}}{\varphi}_{D}(x,y)=0$,
that is, $D$ is $3-$monotone.

Conversely, assume that $D$ is $3-$monotone, i.e., for every $z=(x,y)\in D$,
$\overset{{\scriptscriptstyle \langle2\rangle}}{\varphi}_{D}\le0$.
Since $\iota_{D^{\perp}}\ge0$, in particular, that yields $\iota_{(\Pr_{Y}(D))^{\perp}}(x)\le0$,
namely, $x\in(\Pr_{Y}(D))^{\perp}$. We proved $\Pr_{X}(D)\subset(\Pr_{Y}(D))^{\perp}$. 

In this case $D\subset S\times S^{\perp}=\operatorname*{Graph}(\partial\iota_{S})$
is cyclically monotone, where $S=\Pr_{X}(D)^{\perp\perp}$ is linear
$\sigma(X,Y)-$closed.

(iv) The direct implications are plain due to (i), (iii). Conversely,
assume that $D$ is maximal $\infty-$monotone. In particular, $D$
is $3-$monotone. As previously seen in (iii), $D\subset S\times S^{\perp}$,
where $S=\Pr_{X}(D)^{\perp\perp}$. Since $S\times S^{\perp}$ is
$\infty-$monotone, we infer that $D=S\times S^{\perp}=\operatorname*{Graph}(\partial\iota_{S})=D^{\perp}$
is maximal $\infty-$monotone. But, according to (\ref{eq:-78}) and
Theorem \ref{F}, 
\[
\overset{{\scriptscriptstyle \langle2\rangle}}{\varphi}_{S\times S^{\perp}}(x,y)=\iota_{S}(x)+\iota_{S^{\perp}}(y)=\overset{{\scriptscriptstyle \langle\infty\rangle}}{\varphi}_{S\times S^{\perp}}(x,y),\ (x,y)\in Z.
\]
We conclude that $D$ is maximal $3-$monotone with the aid of Theorem
\ref{fzv}, after we note that $D=[\overset{{\scriptscriptstyle \langle2\rangle}}{\varphi}_{D}=c]$
and $\overset{{\scriptscriptstyle \langle2\rangle}}{\varphi}\ge c$.
\end{proof}

\begin{corollary} Let $(X,Y,\langle\cdot,\cdot\rangle)$ be a dual
system and let $A:D(A)\subset X\rightrightarrows Y$ be a multi-valued
linear skew operator, that is, for every $x\in D(A)$, $y\in A(x)$
, $\langle x,y\rangle=0$. Then

\medskip{}

\noindent\emph{(i)} $A$ is maximal monotone iff $\operatorname*{Graph}A$
is $\sigma(Z,Z)-$closed in $Z$ and $A^{*}:X\rightrightarrows Y$
is monotone, where the adjoint $A^{*}$ is defined by $\operatorname*{Graph}(A^{*})=\neg(\operatorname*{Graph}A)^{\perp}$,
or 
\begin{equation}
(x,y)\in\operatorname*{Graph}(A^{*})\ \Leftrightarrow\ \forall(u,v)\in\operatorname*{Graph}A,\ \langle x,v\rangle=\langle u,y\rangle.\label{eq:-79}
\end{equation}

\noindent\emph{(ii)} $A$ is $3-$monotone iff $A$ is cyclically
monotone iff $R(A)\subset D(A)^{\perp}$. 

\medskip

\noindent\emph{(iii)} $A$ is maximal $3-$monotone iff $A$ is maximal
cyclically monotone iff $D(A)$ is closed convex and $R(A)=D(A)^{\perp}$.

\noindent\end{corollary}

\begin{proof} We use Example \ref{SDC} for $D=\operatorname*{Graph}A$,
for which $D(A)=\Pr_{X}(D)$, $R(A)=\Pr_{Y}(D)$. \end{proof}

\eject

\end{document}